\let\origvec\vec  	
\documentclass[graybox]{svmult}
\let\vec\origvec  	

\usepackage{mathptmx}       
\usepackage{helvet}         
\usepackage{courier}        
\usepackage{type1cm}        

\usepackage{amssymb}
\usepackage{amsmath}
\usepackage{amsfonts}
\usepackage{amscd}

\usepackage{fixltx2e}
\usepackage{pifont}
\usepackage{latexsym}
\usepackage{upgreek}
\usepackage{pifont}
\usepackage{nicefrac} 
\usepackage{cite}
\usepackage{xcolor}
\usepackage{fancybox}
\usepackage{float}
\usepackage{empheq}		
\usepackage{subfig}	
\usepackage{calc}

\usepackage{makeidx}         
\usepackage{graphicx}        
\usepackage{multicol}        
\usepackage[bottom]{footmisc}

\makeindex             
                      
\allowdisplaybreaks[4]
\definecolor{myblue}{rgb}{0,0,0.6}         
\definecolor{gray}{rgb}{0.5,0.5,0.5}

\newcommand*{\der}[2]{\frac{\partial #1}{\partial #2}}
\newcommand{\di}{\,\mathrm{d}}                              
\newcommand{\iin}{\;\text{in}\;}
\newcommand{\oon}{\;\text{on}\;}
\newcommand{\deO}{{\partial\Omega}}
\newcommand*{\conj}[1]{\overline{#1}}

\newcommand*{\N}[1]{\left\|#1\right\|}
\newcommand*{\abs}[1]{\left|#1\right|}
\newcommand{\Tnorm}[1]{|||#1|||}
\newcommand*{\jmp}[1]{[\![#1]\!]}                     
\newcommand*{\mvl}[1]{\{\!\!\{#1\}\!\!\}}             

\DeclareMathOperator{\re}{Re} %
\DeclareMathOperator{\im}{Im} %
\DeclareMathOperator{\diam}{diam} %

\DeclareMathOperator{\essinf}{ess\,inf}

\newcommand{\Uu}[1]{{\mathbf{#1}}}                   

\newcommand{\IC}{\mathbb{C}}

\newcommand{\IN}{\mathbb{N}}
\newcommand{\IR}{\mathbb{R}}

\newcommand{\ba}{{\Uu a}}\newcommand{\bb}{{\Uu b}}
\newcommand{\bd}{{\Uu d}}
\newcommand{\bg}{{\Uu g}}

\newcommand{\bn}{{\Uu n}}

\newcommand{\bu}{{\Uu u}}
\newcommand{\bv}{{\Uu v}}\newcommand{\bw}{{\Uu w}}\newcommand{\bx}{{\Uu x}}
\newcommand{\by}{{\Uu y}}
\newcommand{\bA}{{\Uu A}}
\newcommand{\bF}{{\Uu F}}

\newcommand{\bV}{{\Uu V}}

\newcommand{\bsigma}{{\boldsymbol \sigma}}\newcommand{\btau}{{\boldsymbol\tau}}

           \newcommand{\bzero}{\Uu{0}}

\newcommand{\calA}{{\mathcal A}}\newcommand{\calB}{{\mathcal B}}

\newcommand{\calF}{{\mathcal F}}

\newcommand{\calO}{{\mathcal O}}

\newcommand{\calT}{{\mathcal T}}

\newcommand{\tta}{{\tt a}}
\newcommand{\ttb}{{\tt b}}
\newcommand{\ttd}{{\tt d}} 

\newcommand{\OO}{{(\Omega)}}

\newcommand{\ee}{{\rm e}}
\newcommand{\ri}{{\rm i}}

\newcommand{\beq}{\begin{equation}}      \newcommand{\eeq}{\end{equation}}
\newcommand{\beqs}{\begin{equation*}}    \newcommand{\eeqs}{\end{equation*}}
\newcommand{\bit}{\begin{itemize}}       \newcommand{\eit}{\end{itemize}}
\newcommand{\ben}{\begin{enumerate}}     \newcommand{\een}{\end{enumerate}}
\newcommand{\bal}{\begin{align}}         \newcommand{\eal}{\end{align}}
\newcommand{\bals}{\begin{align*}}       \newcommand{\eals}{\end{align*}}
\newcommand{\bse}{\begin{subequations}}	 \newcommand{\ese}{\end{subequations}}
\newcommand{\bpr}{\begin{proposition}}   \newcommand{\epr}{\end{proposition}}
\newcommand{\bre}{\begin{remark}}        \newcommand{\ere}{\end{remark}}
\newcommand{\bpf}{\begin{proof}}         \newcommand{\epf}{\end{proof}}
\newcommand{\ble}{\begin{lemma}}         \newcommand{\ele}{\end{lemma}}
\newcommand{\bco}{\begin{corollary}}     \newcommand{\eco}{\end{corollary}}
\newcommand{\bex}{\begin{example}}       \newcommand{\eex}{\end{example}}
\newcommand{\bth}{\begin{theorem}}       \newcommand{\enth}{\end{theorem}}
            
\newcommand{\Fh}{\calF_h}
\newcommand{\Th}{{(\calT_h)}}
\newcommand{\dn}{{\partial_\bn}}
\newcommand{\dnK}{{\partial_{\bn_K}}}
\newcommand{\deK}{{\partial K}}
\newcommand{\GD}{{\Gamma_D}}
\newcommand{\GN}{{\Gamma_N}}
\newcommand{\GR}{{\Gamma_R}}

\newcommand{\vhp}{{v_{hp}}}

\DeclareMathOperator*{\argmin}{arg\,min}

\newcommand{\Stdg}{{\textsc{tdg}}}  
\newcommand{\Svtcr}{{\textsc{vtcr}}}  
\newcommand{\Swbm}{{\textsc{wbm}}}  
\newcommand{\Sls}{{\textsc{ls}}}  
  
\newcommand{\Sdgm}{{\textsc{dgm}}}

\newcommand{\inss}{{\mathrm{in}}}
\newcommand{\outss}{{\mathrm{out}}}

\begin{document}

\title*{A Survey of Trefftz Methods for the Helmholtz Equation}
\author{Ralf Hiptmair, Andrea Moiola and Ilaria Perugia}
\institute{Ralf Hiptmair \at Seminar for Applied Mathematics, ETH Z\"urich, 8092 Z\"urich, Switzerland, \email{hiptmair@sam.math.ethz.ch}
\and
Andrea Moiola \at Department of Mathematics and Statistics, University of Reading, Whiteknights PO Box 220, RG6 6AX, UK, \email{a.moiola@reading.ac.uk}
\and
Ilaria Perugia \at Faculty of Mathematics, University of Vienna, 1090 Vienna, Austria, and Department of Mathematics, University of
  Pavia, 27100 Pavia, Italy, \email{ilaria.perugia@univie.ac.at}}
%
%
\maketitle


\abstract{Trefftz methods are finite element-type schemes whose test and trial
    functions are (locally) solutions of the targeted 
    differential equation. They are particularly popular for time-harmonic wave
    problems, as their trial spaces contain oscillating basis functions and
    may achieve better approximation properties than classical
    piecewise-polynomial spaces.\\[\parskip]\hspace*{\parindent}
    We review the construction and properties of several Trefftz variational
    formulations developed for the Helmholtz equation, including least squares,
    discontinuous Galerkin, ultra weak variational formulation, variational theory
    of complex rays and wave based methods. The most common discrete Trefftz
    spaces used for this equation employ generalised harmonic polynomials
    (circular and spherical waves), plane and evanescent waves, fundamental
    solutions and multipoles as basis functions; we describe theoretical and
    computational aspects of these spaces, focusing in particular on their
    approximation properties.\\[\parskip]\hspace*{\parindent}
   One of the most promising, but
    not yet well developed, features of Trefftz methods is the use of adaptivity
    in the choice of the propagation directions for the basis functions.  The main
    difficulties encountered in the implementation are the assembly and the
    ill-conditioning of linear systems, we briefly survey some strategies that
    have been proposed to cope with these problems.}

\section{Introduction}
\label{s:Intro}


Given a linear PDE, a Trefftz method is a volume-oriented discretisation
scheme, for which all trial and test functions, when restricted to any element of a
given mesh, are solutions of the PDE under consideration.  The name comes from the
work \cite{Trefftz1926} of E.~Trefftz, dating back to 1926, where this idea was
applied to the Laplace equation.  Since then, several versions of Trefftz methods
have been proposed and applied to a range of PDEs by different groups of
mathematicians, engineers and computational scientists, often unaware of each
other.  Typical PDEs addressed are \emph{linear}, with \emph{piecewise-constant
  coefficients} and \emph{homogeneous}, i.e.\ with vanishing volume source
term. 

Trefftz methods are related to both finite element (FEM) and boundary element
methods (BEM). With the former they have in common that they provide a
discretisation in the volume. With the latter they share some characteristics such
as the need of integration on lower-dimensional manifolds only. Compared to
conventional FEMs, Trefftz methods have attracted attention mainly for two reasons:
\emph{(i)} they may need much fewer degrees of freedom than standard schemes to
achieve the same accuracy, and \emph{(ii)} they incorporate some properties of the
problem's solution (such as oscillatory character, wavelength, maximum principle,
boundary layers) in the trial spaces, and thus also in the discrete solution.  In
addition, compared to BEMs, an advantage of Trefftz schemes is that they do not
require the evaluation of singular integrals.

Comparing with finite and boundary elements, in 1997 Zienkiewicz \cite{Zie97}
stated: \emph{``\ldots it seems without doubt that in the future Trefftz type
  elements will frequently be encountered in general finite element codes.\ldots
  It is the author's belief that the simple Trefftz approach will in the future
  displace much of the boundary type analysis with singular kernels.''}  While
this prediction has not yet come true, in the last years more and more work has
been devoted to the formulation, the analysis and the validation of these methods
and substantial progress 
has been accomplished.

In this chapter we survey Trefftz finite element methods for the \emph{homogeneous
  Helmholtz equation} ($-\Delta u-k^2 u=0$), which models acoustic wave
propagation in time-harmonic regime.  For medium and high frequencies, i.e.\ for
values of ${k}L$ in a range of $10^{2}$ to $10^{4}$, where $k>0$ is the
wavenumber, and $L$ a characteristic length of the region of interest, the
numerical solution of the Helmholtz equation in 2D and 3D is particularly
challenging.  A main reason is that Helmholtz solutions oscillate with a
wavelength proportional to the inverse of $k$.  Hence, piecewise polynomials do
not provide efficient approximation. Trefftz schemes are thus particularly
relevant as they can improve on the point where (polynomial) FEMs fail: the
approximation properties of the basis functions.  Moreover, some Trefftz methods
can remedy other shortcomings that often haunt discretisations of time-harmonic
problems, such as the lack of coercivity and the presence of minimal resolution
conditions to guarantee unique solvability.  Theorem~\ref{thm:TDG} in this chapter
is an example.  Earlier overviews of Trefftz schemes for the Helmholtz equation,
together with numerous references, can be found in \cite{PHV07},
\cite[Ch.~1]{AndreaPhD} and \cite[Ch.~3]{LuostariPhD}. 
Surveys of Trefftz schemes for other equations are in \cite{Zie97,LLHC08,Qin05,KK95}.

For most of the Trefftz spaces used, continuity across interfaces separating
mesh elements cannot be enforced strongly, as Trefftz functions are not as
``flexible'' as piecewise polynomials.
As a consequence, the standard Helmholtz variational formulation posed in
subspaces of the Sobolev space $H^1$ is not applicable and discretisations
must be used that can accommodate discontinuous trial functions. A wide array of
different variational formulations has been proposed and in \S\ref{s:Methods} we
attempt a classification and a comparison of the best known.  We identify three
main classes of formulations: \emph{(i)} \emph{least squares} (LS, \S\ref{s:LS}),
where squares of suitable norms of residuals are minimised; \emph{(ii)}
\emph{discontinuous Galerkin} (DG, \S\ref{s:DG}), whose formulations arise from
local integration by parts and which may or may not use Lagrange multipliers on
mesh interfaces; \emph{(iii)} \emph{weighted residual} (\S\ref{s:WR}), which are
defined by testing residuals against suitable traces of test functions.  The
methods discussed include: 
the Trefftz-discontinuous Galerkin (TDG), the ultra weak variational formulation
(UWVF), the discontinuous enrichment method (DEM), the variational theory of
complex rays (VTCR) and the wave based method (WBM).
Moreover, in the spirit of the symposium that led up to the present volume, to ``build bridges'' with a wider portion of the literature and of the computational PDE community,
in \S\ref{s:OldStyle} we describe some older Trefftz schemes defined on a single element and in \S\ref{s:PUM_et_al} we consider some methods that are not Trefftz but use oscillating basis functions that are ``approximately Trefftz'', such as the partition of unity method (PUM).
To easily compare them, we write all formulations for the same Robin--Dirichlet
model boundary value problem (see \S \ref{s:BVP_notation}).

In \S\ref{s:Methods} we completely gloss over the choice of \emph{basis functions
  and discrete spaces} employed, whose description is postponed to
\S\ref{s:Spaces}.  This is because, apart from few exceptions such as unbounded
elements, any Trefftz discrete space can be employed in any Trefftz variational
formulation.  We believe that separating the discussion of the two main components
in the definition of a Trefftz method, i.e.\ variational formulations and discrete
spaces, will make the presentation clearer.  The most common basis functions for
Trefftz methods are plane waves ($\bx\mapsto{\rm e}^{\ri k\bd\cdot\bx}$ for a
fixed unit vector $\bd$) and generalised harmonic polynomials (i.e.\ circular/spherical waves, products of circular/spherical harmonics and Bessel functions), for which quite a complete
approximation theory exists, see \S\ref{s:GHP}--\ref{s:PW}.
Other basis functions include fundamental solutions, multipoles, evanescent waves and corner waves.
We note that, since the Helmholtz operator is the sum of a second- and a zero-order term, no non-vanishing piecewise-polynomial Trefftz function is possible.

In this chapter we state a few theorems, none of them is entirely new.
Lemma~\ref{lem:MOW} exemplifies the technique of \cite{MOW99} 
to control the $L^2$ norm of Trefftz functions with mesh-dependent norms containing interface jumps.
If a Trefftz method is well-posed in a suitable skeleton norm, this allows to
control the error in the 
volume; we do this for the LS method in Theorem~\ref{th:LS} and for the TDG method (well-posed by Theorem~\ref{thm:TDG}) in Corollary~\ref{cor:L2TDG}.
This can be combined with the approximation results for circular/spherical and plane waves in \S\ref{s:GHP}--\S\ref{s:PW}. 
In brief: we provide the tools to derive stability and orders of $L^{2}$-convergence in
the 
volume for all Trefftz methods that are well-posed in suitable skeleton norms.

Trefftz methods suffer from two main problems: \emph{ill-conditioning} due to
the poor linear independence of the basis functions, and the need for
\emph{numerical quadrature for oscillating integrands}.  On the other
hand, since the PDE is solved exactly in each element, only low-dimensional
integrals on the mesh skeleton need to be evaluated, leading to massively reduced
computational cost for the assembly of the linear systems.
Moreover, if plane wave bases are used, on
any polygonal/polyhedral mesh the integrals can be computed analytically in a
cheap way. In \S\ref{s:Further} we briefly review strategies developed to deal
with the computation of matrix entries and to cope with ill-conditioning.

Some Trefftz methods also provide an attractive framework for implementing
non-standard adaptive policies, like directional adaptivity following dominant
wave directions. This is made possible, because plane wave-type Trefftz functions
naturally encode a direction of propagation. More details are given in
\S\ref{ss:adapt}.

As mentioned, in this chapter we only discuss 
the Helmholtz equation, i.e.\ acoustic problems, and constant material parameters.
The discrete Trefftz spaces used for the Helmholtz equation with variable
coefficients are briefly addressed 
in \S\ref{s:OTHERS}.  Other time-harmonic
wave problems that have been tackled with Trefftz methods include electromagnetism
(Maxwell equations) \cite{AndreaPhD,CEPhd}, 
linearised Euler equation and general hyperbolic systems \cite{GAB07},
linear elasticity (Navier equation) \cite{LuostariPhD},
(fourth order) Kirchhoff--Love plates \cite{RLK13,LBRK12,DesmetChapter2012,LuostariPhD}, 
Koiter's linear shell theory \cite{RLK13},
poro-elasticity \cite[\S5.4]{DesmetChapter2012}, 
coupled vibro-acoustic problems \cite{DesmetChapter2012}.
A list of applications and references can be found in \cite[\S5.1]{Deckers2014} (with a focus in vibrational mechanics) and in \cite{AndreaPhD,LuostariPhD}.
A related application is tackled by the \emph{method of particular solutions} (MPS) of \cite{FHM67,BeTr05}, which uses Helmholtz solutions to approximate Laplace eigenvalue problems; in this setting the wavenumber is part of the unknowns.
For  recent work on space--time Trefftz methods for wave propagation in time-domain see \cite{SpaceTimeTDG} and references therein.

Several comparisons of the numerical performances of different Trefftz schemes
for simple model problems have been published, e.g.\ \cite{ASG05} (PUM, DEM, generalised FEM), 
\cite{GAA07} (LS, UWVF), \cite{HGA08} (PUM, UWVF), \cite{GGH11} (DG, UWVF, LS),
\cite{WTTF12} (DEM, UWVF, PUM), \cite{HuYuan14} (LS, UWVF, VTCR), where we have included the PUM even if strictly speaking it is not a Trefftz method.  
However, from these results it is difficult to conclude that any
formulation is clearly preferable from a computational point of view.  A general
conclusion might be that, in order to achieve the best accuracy and conditioning,
the choice of the approximation space matters more than that of the variational
formulation.  We reiterate that these two choices are mutually independent: any
Trefftz discrete space might be used in any Trefftz variational formulation.
We make some further concluding remarks in \S\ref{s:concl}.


\subsection{Model boundary value problem}
\label{s:BVP_notation}

We rely on a simple model boundary value problem (BVP) for the Helmholtz equation that will be used to describe and compare the different Trefftz methods. 
Let $\Omega\subset\IR^n$, $n=2,3$, be a bounded, Lipschitz,
connected domain, with $\deO=\GD\cup\GR$, where $\GD$ and $\GR$ are disjoint
components of $\deO$; $\GR\ne\emptyset$ while $\GD$ might be empty.  Denote by
$\bn$ the outward-pointing unit normal vector field on $\deO$.  We consider the
homogeneous Robin--Dirichlet BVP
\begin{equation}\label{eq:BVP}
\begin{aligned}
-\Delta u-k^2 u &=0 && \iin \Omega,\\
u &= g_D	 	&& \oon \GD,\\
\der u\bn + \ri k\vartheta u &= g_R && \oon \GR.
\end{aligned}
\end{equation}
Here $g_D$ and $g_R$ are the boundary data, $\ri$ is the imaginary unit,
$k\in\mathbb{R}$ (the wavenumber) and $\vartheta$ (the impedance parameter) are
positive constants.  We assume that 
$\Omega$, $g_D$ and $g_R$ are such that $u\in H^{3/2+s}\OO$, for some $s>0$.
In typical sound-soft acoustic scattering problems, 
  $\Gamma_D$ represents the boundary of the scatterer, 
  and $\Gamma_R$ stands for an artificial truncation of the unbounded
    region where waves propagate;
see e.g.\ \cite[\S2]{LocallyRef}.

Simple generalisations of the BVP \eqref{eq:BVP} that can be tackled by Trefftz methods are:
\begin{itemize}
\item Neumann boundary conditions $\partial u/\partial \bn=g_N$ on $\Gamma_{D}$; 
\item discontinuous and piecewise-constant wavenumber $k$;
\item piecewise constant and discontinuous tensor coefficient
    $\boldsymbol{A}$
  in the more general Helmholtz
  equation $-\nabla\cdot(\boldsymbol{A}\nabla u)-k^2 u=0$,
    e.g.\ \cite{HMK02} and \cite[Ch.~I.5]{CEPhd};
\item spatially varying impedance $0<\vartheta\in L^\infty(\GR)$;
\item absorbing media $k\in\IC$;
\item inhomogeneous Helmholtz equation $-\Delta u-k^2 u=f$, where the source term $f$ might be either localised \cite[\S5]{GAB07}, \cite{Deckers2014,HowarthPhD,HowarthChildsMoiola2014}, or not 
\cite[\S2.2]{AV05};
\item scattering in unbounded domains; 
\item scattering by periodic diffraction gratings in \cite{ZMZ13,SimonSteve2015Chapter};
\item scattering by screens (i.e.\ manifolds with boundary, leading to non-Lipschitz computational domains) in \cite{ZMZ14}.
\end{itemize}
The presence of smoothly varying coefficients is more challenging for Trefftz methods, as in general no Trefftz functions in analytical form are available; this extension is briefly addressed in \S\ref{s:OTHERS}.

\subsection{Notation}
\label{s:not}

We introduce a finite element partition $\calT_h=\{K\}$ of $\Omega$, not necessarily conforming.
We write $\bn_K$ for the outward-pointing unit normal vector on $\partial K$, and
$h$ for the mesh width of $\calT_{h}$, i.e. $h:=\max_{K\in\calT_h}h_K$, with $h_K:=\diam K$. 
We denote by $\Fh:=\bigcup_{K\in\calT_h}\partial K$ and 
$\Fh^I:=\Fh\setminus\deO$ the skeleton of the mesh and its inner part. 

We also introduce some standard DG notation. 
Given two elements $K_1,K_2\in\calT_h$, a piecewise-smooth function $v$ and vector field $\btau$ on $\calT_h$, we define on $\partial K_1\cap\partial K_2$
\begin{align*}
  \text{the averages:}\quad &
  \mvl{v} := \tfrac12(v_{|K_1} + v_{|K_2}),&&
  \mvl{\btau} := \tfrac12(\btau_{|K_1} +  \btau_{|K_2}),\\
  \text{the normal jumps:}\quad &
  \jmp{v}_N := v_{|K_1}\bn_{K_1}+v_{|K_2}\bn_{K_2},&&
  \jmp{\btau}_N := \btau_{|K_1}\cdot\bn_{K_1}+\btau_{|K_2}\cdot  \bn_{K_2}.
\end{align*}
We denote by $\nabla_h$ the element-wise application of the gradient
$\nabla$, and write $\dn=\bn\cdot\nabla_h$ on $\deO$ and $\dnK=\bn_K\cdot\nabla_h$
on $\partial K$ for the normal derivatives.

For $s>0$, define the broken Sobolev space $H^s\Th$ and the \emph{Trefftz space} $T\Th$:
\begin{align*}
H^s\Th&:=\big\{v\in L^2\OO:\; v_{|K}\in H^s(K) \; \forall K\in\calT_h\big\},\\
T\Th&:=\big\{v\in H^1\Th:\; -\Delta v-k^2v=0 \iin K \text{ and } \dnK v\in L^2(\partial K) \; \forall K\in\calT_h\big\}.
\end{align*}

The discrete Trefftz space $V_p\Th$ is a finite-dimensional subspace of $T\Th$ and can be represented as 
$V_p\Th=\bigoplus_{K\in\calT_h}V_{p_K}(K)$,
where $V_{p_K}(K)$ is a $p_K$-dimensional subspace of $T\Th$ of functions
supported in $K$.  We use the terms $h$-con\-ver\-gence to mean the convergence of
a sequence of numerical solutions to $u$ when the mesh $\calT_h$ is refined, i.e.\
$h\to0$, 
$p$-convergence to designate the convergence when the local spaces are
enriched, i.e.\ $p:=\min_{K\in\calT_h}p_K\to \infty$, and $hp$-convergence to mean
the convergence for a suitable combination of the two refinement strategies.  We
remark that when non-polynomial spaces are used, as it is the case for Trefftz
methods in frequency domain, it is not obvious how to define the ``degree'' of a
space, thus $p_K$ denotes the local number of degrees of freedom.
Finally, we denote by $\re\{\cdot\}$, $\im\{\cdot\}$ and $\conj{\,\cdot\,}$ the
real part, the imaginary part and the conjugate of a complex value.

We note that some of the methods in \S\ref{s:Methods}, such as the TDG, the UWVF
and the VTCR, involve sesquilinear forms (i.e.\ test functions are conjugated)
while others, such as the DEM and the WBM, involve bilinear forms (test functions
are not conjugated).  Any method (if no unbounded elements are used) can be
modified to either form, even though sesquilinear forms are more amenable to
stability and error analysis; for each method we follow the conventions of the
references we cite.

\subsection{Estimation of the $L^2(\Omega)$ norm of (piecewise) Trefftz functions}
\label{s:MOW}

Given two uniformly positive functions $\lambda\in L^\infty(\Fh^I\cup\GD)$ and $\sigma\in L^\infty(\Fh^I\cup\GR)$, we
  introduce the following \emph{skeleton} seminorm (defined e.g.\ on
  $H^{3/2+\varepsilon}\Th$, $\varepsilon>0$):
\begin{align}\label{eq:LambdaSigmaNorm}
\Tnorm{v}^2_{\lambda,\sigma}:=
&\N{\sigma\jmp{\nabla_h v}_N}_{L^2(\Fh^I)}^2
+\N{\lambda\jmp{v}_N}_{L^2(\Fh^I)}^2\\
&+\N{\sigma (\dn v+\ri k\vartheta v)}_{L^2(\GR)}^2
+\N{\lambda v}_{L^2(\GD)}^2.
\nonumber
\end{align}
A special property of the Trefftz space $T\Th$ is that this seminorm is actually a norm for it, and that it controls the $L^2(\Omega)$ norm, as it was first proved by P.~Monk and D.Q.~Wang using a special duality technique in \cite[Th.~3.1]{MOW99}.
\begin{lemma}\label{lem:MOW}
$\Tnorm{\cdot}_{\lambda,\sigma}$ is a norm in $T\Th$.
Moreover, all Trefftz functions $v\in T\Th\cap
H^{3/2+\varepsilon}\Th$, $\varepsilon>0$, satisfy the estimate 
$$\N{v}_{L^2\OO}\le C_*\Tnorm{v}_{\lambda,\sigma},$$
with a constant $C_*>0$ depending only on $k,\lambda,\sigma,\vartheta,\Omega$ and $\calT_h$.
Setting
\[
\sigma_K:=\essinf_{\bx\in\partial
  K\setminus\Gamma_D}\sigma(\bx),\quad 
\lambda_K:=\essinf_{\bx\in\partial
  K\setminus\Gamma_R}\lambda(\bx)\qquad\forall K\in\calT_K,
\]
we can express the dependence of $C_*$ on the relevant parameters in the following situations:

\begin{itemize}
\item[ (i)]\ 
If $\deO=\GR$ and $\Omega$ is either convex or smooth and star-shaped with respect to a ball, then 
$$\N{v}_{L^2\OO}\le 
C_1\,\diam \Omega\,
\max_{K\in\calT_h} \bigg( \Big(\frac1{\sigma_K^2 k}+\frac k{\lambda_K^2}\Big)
\Big(1+\frac1{kh_K}\Big)\bigg)^{1/2}\Tnorm{v}_{\lambda,\sigma},$$
where $C_1>0$ depends on $\vartheta$, the shape-regularity of the mesh and the shape of~$\Omega$.

\item[(ii)]\  
If $k>1$, $\Omega\subset \IR^2$ has diameter $\diam\Omega=1$ and satisfies 
\begin{equation}\label{eq:StarShape}
\bx\cdot\bn\ge\gamma>0 \quad \text{a.e.\ on }\GR  \text{ and}\quad
\bx\cdot\bn\le0\quad \text{a.e.\ on }\GD,
\end{equation}
and each element $K$ is  star-shaped with respect to a ball of radius $\rho_K h_K$, 
we have
\begin{align*}
\N{v}_{L^2\OO}\le C_2 \max_{K\in\calT_h}\bigg(
\Big(\frac{1}{\sigma_K^2 k}+\frac{k}{\lambda_K^2}\Big)
\Big((k h_K)^{2t}+\frac{1}{k h_K}\Big)\bigg)^{1/2}
\Tnorm{v}_{\lambda,\sigma},
\end{align*}
where  
$0<t<s_\Omega\le 1/2$, 
$s_\Omega$ being 
the ``elliptic regularity parameter'' of \cite[eq.~(6)]{LocallyRef},
and $C_2>0$ depends only on $\Omega$, $\vartheta$, $t$, and on the shape-regularity $\inf_{K\in\calT_h}\rho_K$ of the mesh.
\end{itemize}
\end{lemma}
The bound in part \emph{(i)} of Lemma~\ref{lem:MOW} can be verified following the proof of~\cite[Lemma~4.3.7]{AndreaPhD}, while that in part \emph{(ii)} requires also the stability and trace estimates of \cite[eq.~(7), (20)]{HMP13} (see also \cite[Lemma~4.5]{HMP13} and a weaker but more general bound in \cite[Lemma~4.4]{LocallyRef}).
Conditions \eqref{eq:StarShape} on the shape of $\Omega$ are satisfied if $\GR$ is boundary of a domain star-shaped with respect to a ball centred at $\bzero$ and $\GD$ is boundary of a smaller domain (a scatterer, or a ``hole'' in $\Omega$) star-shaped with respect to $\bzero$, 
see \cite[\S2, Fig.~2]{LocallyRef}.
The value of the bounding constants arise only from \emph{(a)} trace estimates for mesh elements, and \emph{(b)} stability bounds for an inhomogeneous Helmholtz BVP on $\Omega$, thus more general shapes of $\Omega$ give different dependencies on $k$ (using e.g.\ the $k$-explicit $H^1\OO$ bounds in \cite[Th.~2.4]{EM11}, \cite[Th.~1.6]{Sp2013a}, and bounds in higher-order norms as in \cite[Lemma~2.12]{GGS14}).
%
This result is relevant because, 
for Trefftz methods that allow a priori stability or error estimates, these are typically in a skeleton norm similar to $\Tnorm{\cdot}_{\lambda,\sigma}$.
Thus Lemma~\ref{lem:MOW} can lead to error estimates in the mesh- and parameter-independent $L^2\OO$ norm; we pursue this in \S\ref{s:LS}, \S\ref{s:TDG}.%

\section{Trefftz variational formulations}
\label{s:Methods}

\subsection{Least squares (LS) methods} 
\label{s:LS}


Least squares methods are perhaps the simplest kind of Trefftz formulations.  They
allow simple error and stability analysis, are easy to implement, lead to
sign-definite Hermitian (or symmetric) linear systems, at the price of a possibly
worse 
conditioning.  A description of Trefftz LS schemes for the Helmholtz equation with
numerous references is given by M.~Stojek in \cite{STO98a}.  The same method is
named \emph{frameless Trefftz elements} in \cite[\S3.6]{Qin05} and \emph{weighted
  variational formulation} (WVF) in \cite{HuYuan14}.  In \cite{MOW99}, Monk and
Wang proposed the following Trefftz LS method for the BVP \eqref{eq:BVP}:
\begin{align}
\nonumber
\text{find}\qquad
u_\Sls&=\argmin_{\vhp\in V_p\Th} J(\vhp;g_R,g_D),\qquad\text{where}\\
J(v;g_R,g_D):&=
\int_{\Fh^I}\Big(\lambda^2\big|\jmp{v}_N\big|^2+\sigma^2\big|\jmp{\nabla_h v}\big|^2\Big)\di S
\label{eq:LS}\\
&\qquad
+\int_\GR\sigma^2\big|\dn v+\ri k \vartheta v-g_R\big|^2\di S
+\int_\GD\lambda^2\big|v-g_D\big|^2\di S,
\nonumber
\end{align}
where $\jmp{\nabla v}:=\nabla_h v_{|K_1}-\nabla_h v_{|K_2}$ on $\partial K_1\cap\partial K_2$ is the jump of the complete gradient (whose ``sign'' depends on a choice of the ordering of the elements in $\Fh$). 
The LS methods in \cite[eq.~(7)]{STO98a} and \cite[Ch.~10]{LLHC08} differ from
\eqref{eq:LS} (apart from the use of different boundary conditions) in that only
the normal component of the jump of the gradient $\jmp{\nabla_h v}_N$ is penalised
on $\Fh^I$, as opposed to the entire jump $\jmp{\nabla_h v}$. 
Obviously, every Galerkin discretisation of the variational problem arising from \eqref{eq:LS} will give rise to a Hermitian linear system, which is a clear advantage of LS methods.

The choice of the relative weights $0<\lambda,\sigma\in L^\infty(\Fh)$ between the terms in \eqref{eq:LS} is a crucial point for the conditioning and the accuracy of LS methods.
Different choices have been proposed (for 2D problems):
$\sigma=1$ and $\lambda=k$ or $\lambda_{|e}=1/h_e$ in \cite[\S2]{MOW99};
$\lambda=1$ and $\sigma_{|e}=h_e/(p_{K_1}+p_{K_2})$ in \cite[\S3.2]{STO98a};
$\lambda=1$ and 
$\sigma_{|e}=\calO(\max\{p_{K_1},p_{K_2}\}^{-1/2})$ in
\cite[Th.~10.3.4]{LLHC08}. 
Here, $e=\partial K_1\cap\partial K_2$ denotes a mesh interface, $h_e$ its length, $p_{K_1}$ and $p_{K_2}$ the dimensions of the local Trefftz spaces $V_{p_{K_1}}(K_1)$ and $V_{p_{K_2}}(K_2)$ on the adjacent elements $K_1$ and $K_2$. 
In 2D and 3D, \cite{HuYuan14} suggests to choose $\sigma=1$ and $\lambda=k$ and, for BVPs with singular solutions, $\sigma_{|\GR}=k^{1/2}$. 

%
%

The LS method 
computes the element $u_\Sls$ in  $V_p\Th$ that minimises the error $u-u_\Sls$
measured in the skeleton
norm $\N{v}_\Sls^2:=J(v;0,0)$, thus orders of converge in this norm follow
immediately from approximation bounds for the specific discrete Trefftz space
$V_p\Th$ chosen, see e.g.\ \S\ref{s:Spaces} below or \cite{MOW99}. 
Since $\Tnorm{v}_{\lambda,\sigma}\le \N{v}_\Sls$ (with equality if $J$ in \eqref{eq:LS} is defined with $\jmp{\nabla_h v}_N$ instead of $\jmp{\nabla_h v}$),
Lemma~\ref{lem:MOW}, following  \cite[Th.~3.1]{MOW99}, guarantees that the $L^2\OO$ norm of the error of the LS solution is controlled by the value of the LS functional, 
thus convergence follows also in $\Omega$.
This is summarised in Theorem~\ref{th:LS}, see \S\ref{s:MOW} for the extension to different domains.
\begin{theorem}\label{th:LS}
Let 
$u$ be the solution of~\eqref{eq:BVP} 
and $u_\Sls\in V_p\Th$ the
discrete LS solution of \eqref{eq:LS}. 
Then, for $C_*>0$ depending only on $k,\lambda,\sigma,\vartheta,\Omega$ and $\calT_h$,
\begin{align*}
\N{u-u_\Sls}_\Sls &=\inf_{\vhp\in V_p\Th}\N{u-\vhp}_\Sls,\\
\N{u-u_\Sls}_{L^2(\Omega)}&\le C_*\, 
\inf_{\vhp\in V_p\Th}\N{u-\vhp}_\Sls.
\end{align*}
If $\lambda=k$, $\sigma=1$, $\deO=\GR$ and $\Omega$ is either convex or smooth and star-shaped,
then
\begin{align*}
&\N{u-u_\Sls}_{L^2(\Omega)}
\le C_0\,
\diam\Omega\,k^{-1/2}\,
\Big(1+\big(k\min_{K\in\calT_h}h_K\big)^{-1/2}\Big)
\inf_{\vhp\in V_p\Th}\N{u-\vhp}_\Sls,
\end{align*}
where $C_0>0$ depends only on $\vartheta$, 
the shape of $\Omega$ 
and the shape-regularity of $\calT_h$.
\end{theorem}
The $hp$-convergence theory of \cite{HMP13} easily extends to the LS method.
  In 2D, if the LS parameters are defined as
  $\lambda^2_{|e}=kh/\min\{h_{K_1},h_{K_2}\}$ for $e=\deK_1\cap \deK_2$,
  $\lambda^2_{|e}=kh/h_{K}$ for $e\subset\deK\cap\GD$, and $\sigma^2=1/k$, under
  the assumptions on $\Omega$ and on the discretisation stipulated in
  \cite{HMP13}, then the $\N{\cdot}_\Sls$ norm of the LS error is estimated as in
  \cite[eq.~(48)]{HMP13} and the $L^2\OO$ norm of the same error converges to zero
  exponentially in the square root of the total number of degrees of freedom
  used.

In \cite[Ch.~10]{LLHC08}, the Trefftz LS scheme is analysed for pure Dirichlet
boundary conditions ($\GR=\emptyset$); the crucial parameter in the analysis is
the relative distance between $k^2$ and the closest Dirichlet eigenvalue of
$-\Delta$.  Error bounds in the broken Sobolev norm $H^1\Th$ are derived.


In the numerical tests in \cite{GGH11} and \cite{GAA07}, the LS method appears to
be slightly less accurate than the UWVF (see \S\ref{s:UWVF} below) and a DG
method, all employed with the same discrete space. 
On the other hand, in the examples in \cite{HuYuan14}, the performance of the LS method is comparable to that of the UWVF and considerably better than that of the VTCR.


\subsubsection{The method of fundamental solutions (MFS)}
\label{s:MFS}

A popular class of LS Trefftz methods is the method of fundamental solutions.
A lucid introduction to the MFS for Helmholtz 
problems, together with numerous references, is in \cite{FKM03}.
The MFS is considered a special case of \emph{source simulation technique} in \cite{Och95}.
The characteristic features of the most common form of the MFS are:
\emph{(i)} the domain is not meshed;
\emph{(ii)} the $N$ basis functions are fundamental solutions ($H^{(1)}_0(k|\bx-\by_\ell|)$ in 2D, $\ell=1,\ldots,N$, where $H^{(1)}_0$ is a Hankel function of the first kind and order zero and $\by_\ell\in\IR^2\setminus\overline{\Omega}$, see \S\ref{s:Hankel});
\emph{(iii)} the minimisation of the $L^2(\deO)$ norm of the error is substituted by the minimisation of the squared error over $M\ge N$ points $\bx_j\in\deO$, $j=1,\ldots,M$.
If $M=N$, the MFS is not an LS method but it simply interpolates the boundary conditions with Trefftz functions.

The same method with plane wave bases is compared to the MFS in \cite{AV05}.  
A variant that is popular in acoustics is the \emph{Helmholtz equation least-squares} (HELS) method, which uses spherical-wave and multipole basis functions, see the recent book \cite{Wu15} and references therein.
LS variants of MFS relying on higher order multipoles in addition to simple Hankel functions have a long history in wave simulations \cite[\S2]{MEH02}.

The locations $\by_\ell$ of the basis singularities are either obtained
numerically together with the coefficients multiplying the basis functions using
non-linear LS solvers \cite[eq.~(7)]{FKM03} (leading to a highly adaptive method),
or can be fixed a priori on a smooth boundary in
$\IR^n\setminus\overline{\Omega}$, e.g.\ using complex analysis techniques (in 2D)
as in \cite{BAB08}, or are determined based on heuristic criteria
\cite[\S3]{MEH02}. 

The MFS with fixed nodes can be interpreted as a discretisation of a compact
  transfer operator related to a single layer potential representation. For this
reason it yields ill-conditioned linear systems; however this does not
rule out efficient computations as demonstrated and analysed in \cite{BAB08}
and in~\cite[\S7]{BAB10}.  According to \cite[p.~766]{FKM03}, the larger the
distance between the nodes and $\Omega$, the more ill-conditioned the linear
system and the more accurate the solution (though this might seem counter-intuitive).

A strength of the MFS is its simplicity of implementation, as no mesh is needed
and all geometric information is contained in only two point sets
$\{\by_\ell\}_{\ell=1}^N\subset\IR^n\setminus\overline\Omega$,
$\{\bx_j\}_{j=1}^M\subset \deO$.  Since fundamental solutions satisfy Sommerfeld
radiation condition, the MFS is often used for scattering problems in unbounded
domains.

In \cite{BAB08}, the convergence of the MFS for Dirichlet problems on a circular
domain is analysed in great detail, and a special design of the curve supporting
the fundamental solutions is proposed for general domains with analytic
boundaries. With this choice, extremely accurate and cheap computations are
possible.

In \cite{BAB10}, Barnett and Betcke present a finite element scheme that couples
the LS formulation of~\cite{STO98a} with the MFS in 2D. They consider the
scattering by sound-soft (non-convex) polygons; the total field is approximated inside an
artificial boundary and the scattered field outside of it.  Singular
Fourier--Bessel basis functions depending on the scatterer's corners (see
\S\ref{s:OTHERS}) are used on all elements adjacent to the scatterer, strongly
enforcing the (homogeneous) Dirichlet boundary conditions; due to this, no terms
on $\deO$ appear in the method formulation. 
Exponential orders of convergence are proved. 
The strong enforcement of boundary conditions may be substituted by an LS approach to deal with more general linear boundary conditions, curved boundaries and transmission problems.

\subsection{Discontinuous Galerkin (DG) methods}
\label{s:DG}

The discontinuous Galerkin (DG) methods constitute a wide class of numerical schemes for the approximation of PDEs, employing discontinuous test and trial functions \cite{ABC01}.
A great number of tools for their design, implementation and error
analysis have been devised, so they are a natural setting for Trefftz methods.
In \cite{HMPS14} we showed that when the interior penalty (IP) method,
one of most common DG schemes, is applied to the Laplace equation, the
use of Trefftz spaces (made of harmonic polynomials) offers better
accuracy than standard spaces also in an $hp$-context.
Similar considerations were made in \cite{LiShu2006} for the
$h$-convergence of the local DG (LDG) method.
To our knowledge, no {\em standard} DG variational formulation (e.g.\ any of those in \cite{ABC01}) has been proposed in the literature to discretise time-harmonic problems with Trefftz basis functions.
Possible reasons for this are that
the error analysis of standard DG schemes requires inverse estimates, which are well-known for polynomial spaces but harder in the Trefftz case (however, see \cite[\S3.2]{GHP09} for $h$-explicit inverse estimates for plane waves in 2D),
and that the application of formulations designed for the Laplace equation to the Helmholtz case requires some problematic minimal resolution condition to ensure unique solvability \cite{MPS12}.

In the next subsections we outline some DG formulations that have been
designed specifically for Trefftz discretisations; some of these have
later been employed also with polynomial approximating spaces, e.g.\ \cite{MPS12,MSS10}.

A note on terminology: all Trefftz methods presented in this survey involve the
discretisation of variational formulations based on discontinuous functions,
however with ``DG'' we denote only those methods that arrive at local variational formulations by applying integration by parts to the PDE to be approximated. 
On the contrary, least squares and weighted residual methods simply enforce (weakly) continuity and boundary conditions, irrespectively of the considered PDE.

\subsubsection{The Trefftz-DG (TDG) method}
\label{s:TDG}

Originally, Trefftz-discontinuous Galerkin (TDG) methods (or plane
wave DG, PWDG, when used in combination with plane wave basis functions) were introduced as a way of recasting the ultra weak variational formulation (UWVF) 
of~\cite{CED98,CEPhd} (see \S\ref{s:UWVF} below) in a framework that would facilitate its theoretical analysis \cite{BUM07,GHP09}. 
A similar, but more general, Trefftz-DG framework was proposed in \cite{GAB07,GGH11}, arising from methods for hyperbolic equations; see Remark~\ref{rem:Hyperbolic} below.

We first derive the TDG formulation as in \cite{LocallyRef}.
We multiply the Helmholtz equation \eqref{eq:BVP} by a test function $v$ and integrate by
parts twice on each $K\in\calT_h$:
\begin{align*}
0&=\int_K(-\Delta u-k^2u)\conj v\di V      \quad=\quad
\int_K(\nabla u\cdot\conj{\nabla v}-k^2 u\conj{v})\di V
-\int_{\partial K}\nabla u\cdot\bn_K\,\conj{v}\di S\\
&=\int_K u\,(-\Delta\conj{v}-k^2 \conj{v})\di V+
\int_{\partial K} u\,\conj{\dnK v}\di S
-\int_{\partial K} \dnK u\,\conj{v}\di S.
\end{align*}
We then replace $u$ and $v$ by discrete functions $u_{hp},v_{hp}\in V_p(\calT_h)$, the trace of $u$ on $\partial K$ by the  numerical flux $\widehat{u}_{hp}$, 
and the trace of $\nabla u$ by the numerical flux $\ri k\widehat{\bsigma}_{hp}$ (both defined below),
obtaining the elemental TDG formulation:
\begin{align}\label{eq:DGelemental}
\int_{\partial K}\widehat{u}_{hp}\,\conj{\dnK v}_{hp}\di S
-\int_{\partial K} \ri k\widehat{\bsigma}_{hp}\cdot\bn_K\,\conj{v}_{hp}\di S
=0,
\end{align}
where the volume integral vanishes as the test function $v_{hp}\in V_P\Th\subset T\Th$ is a Trefftz function.
%
%
%
Variants of DG methods are distinguished by the underlying numerical fluxes. Here we opt for the \emph{primal fluxes}:
\begin{align}
\ri k\widehat{\bsigma}_{hp}&=
\begin{cases}
\displaystyle{\mvl{\nabla_h u_{hp}}-\alpha\, \ri k\,\jmp{u_{hp}}_N}&
\text{on faces in $\calF_h^I$},\\
\displaystyle{\nabla_h u_{hp}
-(1-\delta)\left(\nabla_h u_{hp}
+\ri k\vartheta u_{hp}\bn-\,g_R\bn\right)}\hspace{5mm}&
\text{on faces in $\GR$},\\
\displaystyle{\nabla_h u_{hp}-\alpha\,\ri k\,( u_{hp}-g_D)\bn}&
\text{on faces in $\GD$},
\end{cases}\label{eq:TDGfluxes_sigma}
\\
\widehat{u}_{hp}&=
\begin{cases}
\displaystyle{\mvl{u_{hp}}-\beta\,(\ri k)^{-1}\jmp{\nabla_h u_{hp}}_N}&
\text{on faces in $\calF_h^I$},\\
\displaystyle{u_{hp}-\delta\left((\ri k\vartheta)^{-1}\nabla_h u_{hp}\cdot\bn
+u_{hp}-(\ri k\vartheta)^{-1}g_R\right)}&
\text{on faces in $\GR$},\\
g_D &\text{on faces in $\GD$},
\end{cases}\label{eq:TDGfluxes_u}
\end{align}
where the flux  parameters $\alpha>0$, $\beta>0$, $0<\delta\le1/2$, are bounded functions defined on suitable unions of edges/faces (see also Table~\ref{tab:TDGfluxParam}). 
Adding over all elements, we obtain the following formulation of the TDG method: 
\begin{align}
&\text{find } u_\Stdg\in V_p(\calT_h)\text{ s.t. }\quad
\calA_\Stdg(u_\Stdg,v_{hp})=\ell_\Stdg(v_{hp})\quad\forall v_{hp}\in V_p(\calT_h),
\quad\text{where}
\nonumber\\
&\calA_\Stdg(u,v):= 
\label{eq:TDG}\\
&\int_{\calF_h^I}
\Big(\mvl{u}\jmp{\conj{\nabla_h v}}_N
-\mvl{\nabla_h u}\cdot\jmp{\conj{v}}_N
+\alpha \ri k\jmp{u}_N\cdot\jmp{\conj{v}}_N
-\beta(\ri k)^{-1}\jmp{\nabla_h u}_N\jmp{\conj{\nabla_h v}}_N
\Big)\di S
\nonumber
\\
&+\int_\GR \Big(
(1-\delta)\ri k\vartheta u\conj{v}
+(1-\delta)u \conj{\dn v} 
-\delta\dn u \:\conj{v}
-\delta(\ri k\vartheta)^{-1}\dn u \conj{\dn v}
\Big)\di S
\nonumber
\\
&+\int_\GD \Big(-\dn u\:\conj{v} +
\alpha\,\ri k\,u\,\conj{v}\Big)\di S,
\nonumber
\\
&\ell_\Stdg(v):=
\int_\GR g_R\Big((1-\delta)\conj{v}-\delta(\ri k\vartheta)^{-1}
\conj{\dn v}
\Big)\di S
+\int_\GD g_D\Big(\alpha \ri k\conj v-\conj{\dn v}\Big)\di S. 
\nonumber
\end{align}
The TDG method was introduced in the primal form described here in \cite{GHP09,GIT08} and in mixed form in \cite{HIP08p}, under the name of \emph{plane wave DG (PWDG) method}, following the derivation of \cite{ABC01} of general DG schemes for elliptic equations.
In \cite{GHP09}, first-order convergence in the meshwidth was
established, using Schatz' argument, for 2D Robin problems with
  source term $f\in L^2\OO$, plane wave discrete spaces and quasi-uniform
  families of meshes.
This was extended to higher orders in $h$ in \cite{MOI09}, $p$-convergence in \cite{PVersion}, three dimensions in \cite{AndreaPhD}, locally-refined meshes in \cite{LocallyRef}, and finally the exponential convergence in the number of degrees of freedom of its $hp$-version was proved in \cite{HMP13}.
Its dispersion analysis was performed in \cite{GIT08,GitHip14}.

For polynomial discrete spaces, the advantages of using the 
formulation underlying the TDG method, compared to standard DG schemes, were analysed in \cite{MPS12}.
In \cite{BeP12}, the TDG formulation was utilised with (non-Trefftz) bases defined from oscillating functions from high-frequency asymptotics modulated with polynomials; problems with varying coefficients were also considered.

The TDG formulation \eqref{eq:TDG} 
can be seen as a modification of either
the \emph{interior penalty} method, or of the \emph{local DG} (LDG)
method (see e.g.\ \cite{ABC01}):
with respect to the \emph{interior penalty} method, the stabilisation term
multiplied by $\beta$ is added in the TDG fluxes
\eqref{eq:TDGfluxes_u}, while with respect to the LDG method, 
in the TDG fluxes \eqref{eq:TDGfluxes_sigma},
the consistency term is written
in terms of the primal variable ($\mvl{\nabla_h u_{hp}}$) instead of
in terms of the auxiliary variable ($\mvl{ik{\bsigma}_{hp}}$) and 
the additional stabilisation of the jumps of ${\bsigma}_{hp}$ is removed.
In \cite{Spence2015else}, the TDG and the UWVF are seen as special instances of a family of methods arising from integration by parts.

The a priori error analysis of the TDG relies on Theorem~\ref{thm:TDG} below (e.g.\ \cite[\S4]{LocallyRef}), which makes use of the following mesh- and flux-dependent seminorms:
\begin{align*}
\Tnorm{v}_\Stdg^2 := & \;k^{-1}\N{\beta^{\frac12}\jmp{\nabla_h v}_N}_{L^2(\calF_h^I)}^2
+k\N{\alpha^{\frac12}\jmp{v}_N}_{L^2(\calF_h^I)}^2\\
&+k^{-1}\N{\delta^{\frac12}\vartheta^{-\frac12}\dn v}_{L^2(\GR)}^2
+k\N{(1-\delta)^{\frac12}\vartheta^{\frac12}v}_{L^2(\GR)}^2
+k\N{\alpha^{\frac12}v}_{L^2(\GD)}^2;
\\
\Tnorm{v}_{\Stdg^+}^2:=&\;\Tnorm{v}_\Stdg^2+k\N{\beta^{-\frac12}\mvl{v}}_{L^2(\calF_h^I)}^2
+k^{-1}\N{\alpha^{-\frac12}\mvl{\nabla_h v}}_{L^2(\calF_h^I)}^2\\
&+k\N{\delta^{-\frac12}\vartheta^{\frac12}v}_{L^2(\GR)}^2
+k^{-1}\N{\alpha^{-\frac12}\dn v}_{L^2(\GD)}^2.
\end{align*}
\begin{theorem}\label{thm:TDG}
The seminorms $\Tnorm{\cdot}_\Stdg$ and $\Tnorm{\cdot}_{\Stdg^+}$ are norms in the Trefftz space $T\Th$.
The TDG sesquilinear form is continuous and coercive:
\begin{equation*}
\abs{\calA_\Stdg(v,w)}\le 2\Tnorm{v}_{\Stdg^+}\Tnorm{w}_\Stdg,\quad  
\im\big\{\calA_\Stdg(v,v)\big\}=\Tnorm{v}_\Stdg^2 
\end{equation*}
for all $v,w\in T\Th$,
thus there exists a unique solution $u_\Stdg\in V_p\Th$ to the TDG formulation \eqref{eq:TDG} and the quasi-optimality bound holds:
\begin{equation*}
\Tnorm{u-u_\Stdg}_\Stdg\le 3\inf_{v_{hp}\in V_p(\calT_h)}\Tnorm{u-v_{hp}}_{\Stdg^+}.
\end{equation*}
\end{theorem}
Choosing 
$\lambda^2=\alpha k$ on $\Fh^I\cup\GD$, $\sigma^2=\beta/k$ on $\Fh^I$ and $\sigma^2=\min\{\delta,1-\delta\}/2k\vartheta$ on $\GR$,
the norm \eqref{eq:LambdaSigmaNorm} is controlled as
$\Tnorm{v}_{\lambda,\sigma}\le \Tnorm{v}_\Stdg$ for all $v\in T\Th$.  Thus, by
Lemma~\ref{lem:MOW}, the $L^2\OO$ norm of the TDG error can be controlled by its
$\Tnorm{\cdot}_\Stdg$ norm, and so by the discrete space approximation properties.
This result has been stated 
in several slightly different forms, depending on the
regularity of the solution $u$, the type of mesh used, the choice of the numerical
flux parameters $\alpha,\beta,\delta$; see \cite[Lemma~4.3.7]{AndreaPhD},
\cite[Lemma~4.4]{LocallyRef} and \cite[Lemma~4.5]{HMP13}.  
To strike a balance between the size of the constants arising from the duality argument of Lemma~\ref{lem:MOW} and approximation errors, different flux parameters have been chosen
on different meshes and aiming at different types of
convergence estimates, see Table~\ref{tab:TDGfluxParam}.  For illustration, we
state the result in the case of constant flux parameters, quasi-uniform meshes,
and domains
that guarantee sufficiently smooth solutions for the dual problems;
this follows from Lemma~\ref{lem:MOW} and Theorem~\ref{thm:TDG} (\emph{cf.}~\cite[Cor.~4.3.8]{AndreaPhD}).
\begin{corollary}  \label{cor:L2TDG}
Let $u$ be the solution of~\eqref{eq:BVP}, where $\Omega$ is either convex or smooth and star-shaped,
and let $u_{\Stdg}\in V_p(\calT_h)$ be the solution of the TDG method with flux parameters chosen as in the second row of Table~\ref{tab:TDGfluxParam}. 
Then
\begin{align*}
&\N{u-u_{\Stdg}}_{L^{2}(\Omega{})} 
\le C_0 \, \diam\Omega\,   \Big(1+\big(k\min_{K\in\calT_h}h_K\big)^{-1/2}\Big)
\inf_{v_{hp}\in V_p(\calT_h)}\Tnorm{u-v_{hp}}_{\Stdg^+},
\end{align*}
where $C_0>0$ depends only on $\vartheta$, 
the shape of $\Omega$ 
and the shape-regularity of the mesh, but is independent of $k$ and $V_{p}(\mathcal{T}_{h})$.
\end{corollary}


The combination of the abstract error analysis outlined above and approximation estimates for plane, circular and spherical waves (see \S\ref{s:Spaces}) leads to a priori $h$-, $p$- and $hp$-convergence estimates in $\Tnorm{\cdot}_\Stdg$ and $L^2$ norms, see \cite{GHP09,PVersion,AndreaPhD,LocallyRef,HMP13}.
The dependence of the error bounds on the wavenumber $k$ is explicit, as in Corollary~\ref{cor:L2TDG}.

\begin{table}[htb]\centering
\begin{tabular}{|ll|l|l|l|}\hline
&&$\alpha$&$\beta$&$\delta$\\\hline
Quasi-uniform meshes, $h$-convergence & \cite{GHP09} & $\tta/kh_K $ &
$\ttb kh_K $ & $\ttd kh_K $\\
Quasi-uniform meshes, $p$-convergence & \cite{PVersion} & $\tta $ & $\ttb $ & $\ttd $\\
UWVF (see \S\ref{s:UWVF}) & \cite{CED98}&$1/2$&$1/2$&$1/2$\\
Locally refined meshes, $hp$-convergence&\cite{LocallyRef} & $\tta h/h_K $ & $\ttb h/h_K $ & $\ttd h/h_K $\\
Geometrically graded meshes, exponential $hp$-convergence &\cite{HMP13} & $\tta h/h_K $ & $\ttb $ & $\ttd $\\
Polynomial (non Trefftz) basis, $hp$-convergence &\cite{MPS12} & $\tta q_K^2/kh_K $ & $\ttb kh_K/q_K $ & $\ttd kh_K/q_K $
\\\hline
\end{tabular}
\caption{Different TDG flux parameters in \eqref{eq:TDGfluxes_sigma} and
  \eqref{eq:TDGfluxes_u} that have been considered.
  Here $\tta,\ttb,\ttd$ are positive functions independent of the other parameters;
  $k$ is the wavenumber;
  $h_K$ is the local meshwidth;
  $h=\max_{K\in\calT_h}h_K$ is the global meshwidth;
  $q_K$ is the local polynomial degree (for the non-Trefftz version).
}
\label{tab:TDGfluxParam}
\end{table}

\begin{remark}\label{rem:Hyperbolic}
The Helmholtz equation may be written as the first order hyperbolic system
$-\ri k\bu+\sum_{j=1}^n\partial_{x_j}(\bA^{(j)}\bu)=\bzero$, where $\bu:=(u;\nabla u/(\ri k))$ and $\bA^{(j)}$ are the $(1+n)\times(1+n)$ symmetric matrices whose only non-zero elements are $A^{(j)}_{1,j+1}=A^{(j)}_{j+1,1}=1$, for $1\le j\le n$.
Then, similarly to \cite[eq.~(22)]{GAB07} or \cite[eq.~(5)]{GGH11}, a general Trefftz-DG method can be written as:
\begin{align*}
&\quad\text{seek }\;\bu\in \bV_p\Th:=\big\{(u,\bsigma):\; u\in V_p\Th, \bsigma= \nabla u/(\ri k)\big\}\;
\text{ s.t. }\; \forall \bv\in\bV_p\Th
\\
&
\sum_{\substack{K_1,K_2\in\calT_h,\\K_1\ne K_2}}
\int_{\deK_1\cap\deK_2} \!
\big(\bF_{|K_1}^\inss\bu_{|K_1}-\bF_{|K_2}^\inss\bu_{|K_2}\big)\cdot
\big(\conj {\bv_{|K_1}}- \conj{\bv_{|K_2}}\big)\di S
+\!\int_\deO (\bF^\inss \bu-\bg)\cdot\conj \bv\di S=0
\end{align*}
where the flux-splitting matrices $\bF^\inss,\bF^\outss$ are defined on $\prod_{K\in\calT_h}\deK$ and satisfy 
$\bF^\inss\le0$, $\bF^\outss\ge0$ (i.e.\ are negative and positive semi-definite, respectively),
$\bF^\inss+\bF^\outss=(\begin{smallmatrix}0&\bn_K^\top\\ \bn_K &\bzero\end{smallmatrix})$ on $\deK$,
and $\bF_{K_1}^\inss=-\bF^\outss_{K_2}$ on $\deK_1\cap\deK_2$.
The boundary data are represented by a suitable 
vector field $\bg=-\bF^\outss\bu$.
The TDG in \eqref{eq:TDG} (up to a factor $-\ri k$) is obtained by choosing:
\begin{align*}
\begin{aligned}
&\bF^\inss_K=\hspace{37mm}\bF^\outss_K=
\\
&\begin{cases}
\begin{pmatrix}-\alpha&\frac12\bn_K^\top\\ \frac12\bn_K &-\beta\bn\otimes\bn^\top\end{pmatrix}\\
\begin{pmatrix}-(1-\delta)\vartheta&\delta\bn_K^\top\\ (1-\delta)\bn&-\frac\delta\vartheta\bn\otimes\bn^\top\end{pmatrix}\\
\begin{pmatrix}-\alpha&\bn_K^\top\\ \bzero&\bzero\end{pmatrix}
\end{cases}
\begin{cases}
\begin{pmatrix}\alpha&\frac12\bn_K^\top\\ \frac12\bn_K &\beta\bn\otimes\bn^\top\end{pmatrix}
&\oon \deK\cap\Fh^I,\\
\begin{pmatrix}(1-\delta)\vartheta&(1-\delta)\bn_K^\top\\
\delta\bn&\frac\delta\vartheta\bn\otimes\bn^\top\end{pmatrix}
&\oon \deK\cap\GR,\\
\begin{pmatrix}\alpha&\bzero^\top\\ \bn_K &\bzero\end{pmatrix}
&\oon \deK\cap\GD.
\end{cases}
\end{aligned}
\end{align*}
The right-hand side is represented by the vector
$\bg=-\frac1{\ri k}(\begin{smallmatrix}1-\delta\\ \delta\vartheta^{-1}\bn_K \end{smallmatrix})g_R$
on $\GR$
and $\bg=-(\begin{smallmatrix}\alpha\\ \bn_K \end{smallmatrix})g_D$ on $\GD$.
\end{remark}


\subsubsection{The ultra weak variational formulation (UWVF)}
\label{s:UWVF}


The ultra weak variational formulation (UWVF) has been introduced in the 1990's by O.~Cessenat and B.~Despr\'es in~\cite{CED98,CEPhd}.
Since then it has received a great deal of attention and has been applied to numerous PDEs and BVPs; we refer to \cite{HMK02} for a description of its computational aspects and to \cite[\S3.5.2]{LuostariPhD} for an extensive bibliography.
Different derivations can be found e.g.\ in \cite{CED98,BUM07,GHP09,GAB07,GGH11}; 
in particular \cite{GHP09,BUM07} obtain the UWVF in the setting of DG schemes for elliptic problems of \cite{ABC01}, while \cite{GAB07,GGH11}  
derive it for general first-order hyperbolic systems using a flux-splitting approach as we did for the TDG in Remark~\ref{rem:Hyperbolic}.
Note that different papers use different sign conventions.
The extension of the UWVF to problems with smooth coefficients has been tackled in \cite{ImbertGeradDespres2014}.

To write its formulation for the BVP \eqref{eq:BVP} in the Robin case, i.e. $\GD=\emptyset$,
we first define the trace space $X:=\prod_{K\in\calT_h} L^2(\partial K)$, 
and the operators $F_K: L^2(\partial K)\rightarrow  L^2(\partial K)$, mapping the boundary datum  $y_K$ of a local adjoint-impedance Helmholtz BVP into 
the impedance trace of the BVP solution $e_K$ itself:
\begin{align*}
F_K(y_K):=(\partial_{\bn_K}+\ri k)e_K, \qquad\text{where}\qquad
 \begin{cases}
-\Delta e_K-k^2 e_K =0   & \iin K,\\
(-\partial_{\bn_K} +\ri k) e_K = y_{K}  &\oon \partial K.
 \end{cases}
\end{align*}
The Helmholtz BVP is written as a transmission problem across the mesh interfaces, i.e., for all $K,K'\in \calT_h$,
\[
\begin{aligned}
-\Delta u-k^2 u &=0 && \iin K,\\
\partial_{\bn_K}u+\ri k u &= -\partial_{\bn_{K'}}u+\ri k u	 	&& \oon \partial K\cap\partial K^\prime,\\
\partial_{\bn_K}u+\ri k \vartheta u &= g_R && \oon \partial K\cap\GR.
\end{aligned} 
\]
Then, after multiplying the first equation by $e_{|K}$, $e\in T\Th$, integrating by parts twice, taking into account transmission and boundary conditions, and introducing $x,y\in X$ defined as $x_{|\deK}=-\dnK u+\ri k u$ and $y_{|\deK}=-\dnK e+\ri k e$, the
UWVF of problem~\eqref{eq:BVP} \cite[(1.4)]{CED98} reads: find $x\in X$ 
such that, for every $y\in X$,
\begin{align}
\label{eq:UWVF}
&\sum_{K\in\calT_h}\int_{\partial K} x_{|\partial K} \:\conj{y_{|\partial K}}\di S
-\sum_{K,K'\in\calT_h}\int_{\partial K\cap\partial K'} x_{|\partial K'}\: \conj{F_K(y_{|\partial K})}\di S\\
&\qquad-\sum_{K\in\calT_h}\int_{\partial K\cap\GR}
\frac{1-\vartheta}{1+\vartheta} x_{|\partial K}\: \conj{F_K(y_{|\partial K})}\di S
=\sum_{K\in\calT_h}\int_{\partial K\cap\GR} \frac2{1+\vartheta}g_R\: \conj{F_K(y_{|\partial K})}\di S.
\nonumber
\end{align}
(Note that for $\vartheta=1$ the term on $\deK\cap\GR$ at left-hand side vanishes and $2/(1+\vartheta)=1$.)
The expression~\eqref{eq:UWVF} is a variational formulation for the skeleton unknown $x$; after the equation is solved for $x$, the Helmholtz solution $u_{|K}$ 
can be recovered in the interior of each element by solving a local (in $K$) adjoint-impedance Helmholtz BVP with datum $(-\partial_{\bn_K}+\ri k) u_{|K}=x_{|\partial K}$.
If the formulation is discretised choosing a finite dimensional subspace $X_h$ of $X$ corresponding to the impedance traces of a Trefftz space, namely
$$X_h:=\big\{x_h\in X:\, 
{x_h}_{|\partial K}=(-\partial_{\bn_K}+\ri k) v_{|K}\;\forall K\in\calT_h,\; v\in V_p\Th
\big\},$$
then the action of  
$F_K$ and the reconstruction of $u_K$ in $K$ are immediately computed.

Theorem 2.1 of~\cite{CED98} states that the discrete problem obtained by substituting $X_h$ to $X$ in~\eqref{eq:UWVF} is 
solvable, independently of the meshsize $h$;
Corollary~3.8 shows that, for plane wave discrete spaces, the Dirichlet and Robin traces of the UWVF solution converge to the corresponding traces of $u$ with algebraic orders of convergence in $L^2(\GR)$.
In \cite[\S4]{BUM07}, these results have been used together with the duality technique of~\cite{MOW99} to prove orders of convergence for the $L^2\OO$
norm of the error. 

The UWVF has been recast as a DG method with Trefftz basis functions in several
different ways in \cite{BUM07,GAB07,GGH11,GHP09}.  
In particular, \cite[Remark~2.1]{GHP09} shows that \emph{the UWVF is a
    special case of the TDG} formulation \eqref{eq:TDG} for flux parameters
$\alpha=\beta=\delta=1/2$.  As a consequence, the orders of convergence in $h$
and $p$ proved for the TDG on quasi-uniform meshes in \cite{GHP09,PVersion}
carry over to the UWVF (with suboptimal orders in $h$); on the other hand, the
$hp$-type results of \cite{LocallyRef,HMP13} require variable numerical flux
parameters to cope with elements of different sizes (see
Table~\ref{tab:TDGfluxParam}), so they do not apply to the UWVF.
Thus, the TDG can be understood as the extension of the UWVF to non quasi-uniform meshes.
Alternatively, in \cite[\S4.3, 5.2]{MSS10}, the UWVF is employed on meshes refined towards solution singularities by choosing Trefftz spaces on large elements and polynomial spaces on small ones.
No applications of the TDG combining mesh-dependent parameters and polynomial spaces in small elements have been documented.

\subsubsection{DG schemes with Lagrange multipliers}
\label{s:Farhat}

The DG schemes described so far enforce weak continuity between elements using numerical fluxes, in the spirit of \cite{ABC01}.
A different approach is to enforce continuity using Lagrange multipliers.
This was probably first proposed for Trefftz methods  in \cite[\S2.3]{IhBa97}, for the 1D Helmholtz equation.

This strategy has been followed in the \emph{discontinuous enrichment method} (DEM), introduced by C.~Farhat, I.~Harari and L.P.~Franca in \cite{FHF01}, combining a space of piecewise-constant Lagrange
multipliers on mesh interfaces with a discrete space composed by sums of continuous piecewise polynomials and discontinuous plane waves.
Subsequently, in \cite{FHH03}, the polynomial part of the trial space was dropped,
leaving a plane wave trial space and thus
reducing to a Trefftz method; in this version, the DEM was renamed \emph{discontinuous
Galerkin method} (DGM) and the Lagrange multipliers were approximated by
oscillatory functions.  This formulation performed very well for test cases and
was later extended to ``higher order elements'' (i.e.\ elements containing more
plane waves) and other PDEs.  We refer again to \cite[\S3.5.3]{LuostariPhD} for a
comprehensive bibliography.


Here we briefly describe the formulation of the DGM following \cite[\S2]{FHH03}:
\begin{align*}
&\text{find }(u,\lambda)\in H^1\Th\times W\Th \text{ s.t. }\\
&\begin{cases}
\displaystyle\calA_\Sdgm(u,v)+\calB_\Sdgm(\lambda,v)=\int_\GR g_R\,v\di S
\qquad&\forall v\in  H^1\Th,\\
\displaystyle
\calB_\Sdgm(\mu,u)= \int_\GD \mu\,g_D\di S
&\forall\mu\in W\Th,
\end{cases}
\\
&\text{where}
\\
\calA_\Sdgm(w,v):&=\sum_{K\in\calT_h} \int_{K}(\nabla w\cdot\nabla v-k^2 u\,v)\di V + \int_{\GR}\ri k\vartheta \,w\,v\di S,\\
\calB_\Sdgm(\mu,w):&=\sum_{K,K'\in\calT_h} \int_{\partial K\cap\partial K'}\mu(w_{|K'}-w_{|K})\di S
+\int_\GD \mu\,w\di S,
\\
W\Th:&=\bigg(\prod_{K,K'\in\calT_h} 
{\tilde{H}^{-1/2}}(\partial K\cap\partial K')\bigg)\times H^{-1/2}(\GD).
\end{align*}
It is immediate to verify that the solution $u$ to BVP \eqref{eq:BVP} satisfies this formulation, and that the multiplier $\lambda$ represents the normal derivative of $u$ on the mesh interfaces and on $\GD$.
This formulation is then discretised by restricting it to finite dimensional spaces $V_p\Th\subset H^1\Th$ and $W_p\Th\subset W\Th$.
In the DEM of \cite{FHF01}, $V_p\Th$ is the direct sum of a continuous polynomial and a plane wave space, in the DGM of \cite{FHH03} and subsequent papers only the plane wave part is retained, so  $V_p\Th\subset T\Th$.
The volume degrees of freedom, i.e.\ those corresponding to $V_p\Th$, are then eliminated by static condensation in order to reduce the computational cost of the scheme.


A stability and convergence analysis of the simplest version of the DGM (four plane waves per
element and piecewise-constant multipliers) is attempted in~\cite{ADF09}: for
a Robin--Neumann BVP on a domain decomposed in rectangles, under a mesh resolution
condition, the scheme is shown to be well-posed, and a priori orders of
convergence are proved (in $H^1\Th$ norm for the primal variable and in $L^2(\Fh)$
for the multipliers), along with residual-type a posteriori error bounds.
We are not aware of any error analysis for the DGM method holding in more general situations (e.g.\ more than four plane waves per elements, propagation directions not aligned to the mesh, non-rectangular mesh elements).

A similar formulation, named \emph{hybrid-Trefftz finite element method}, is
described in \cite[\S3.5]{Qin05} (deriving the functional in eq.~(65) therein):
the same form $\calA_\Sdgm$ above is used, while $\calB_\Sdgm$ is substituted by
$\calB_{\textsc{ht}}(\mu,w):=-\int_{\Fh^I}\mu\,\jmp{\nabla_h w}_N\di
S-\int_\GN\mu\,\dn w\di S$, where now the multiplier $\mu$ approximates the
Dirichlet trace of $u$, the right-hand sides and the space $W\Th$ are changed
accordingly.  
A further variant of hybrid-Trefftz methods is presented in \cite{SLF10} and related papers.

Another DG method with Trefftz basis, called \emph{modified DG method} (mDGM), has been proposed in \cite{GACD12}.  
The Lagrange multipliers are double-valued on the interfaces (differently from the DEM/DGM of \cite{FHF01,FHH03}) and belong to $\prod_{K\in\calT_h}L^2(\deK\setminus\GR)$.  
A two-step procedure is adopted. 
First, for each basis element $\lambda\in L^2(\deK\setminus\GR)$ of the discrete Lagrange multiplier space, a well-posed Helmholtz BVP on $K$ with impedance datum $\lambda$ is solved in the local Trefftz space $V_{p_K}(K)$ using the classical $H^1(K)$-conforming variational formulation.
Second, these local solutions are combined in a global LS formulation leading to a positive semi-definite system whose unknowns are the Lagrange multipliers themselves.
%
The mDGM was further improved in \cite{ACDG12} leading to the \emph{stable DG method} (SDGM), which differs from the mDGM in that the local impedance problems are solved with a least squares formulation posed on $\deK$, which gives local Hermitian matrices.


Lagrange multipliers are also used to tackle problems with discontinuous coefficients by means of the partition of unity method, see \cite{LBPT05} and \S\ref{s:PUM_et_al} below.

\subsection{Weighted residual methods}\label{s:WR}

Trefftz discretisations lend themselves well to weighted residual formulations: the discrete solution is automatically a local solution of the PDE, only the residual on interfaces (the jumps) and on the boundary (the mismatch with boundary conditions) need to be enforced by multiplying them to suitable traces of test functions.
The choice of these traces leads to different variational formulations, the most developed of which are the VTCR and the WBM described in the following.
While it is simple to design weighted residual methods, their error analysis is by
no means easy, as they arise neither from integration by parts, nor from a minimisation principle.

An earlier weighted-residual Trefftz formulation is the \emph{weak element method} of \cite{Gol86}, where the integral averages of Dirichlet and Neumann jumps on mesh faces are set to zero (equivalently, test functions are constant on each mesh face).

We note that some of the earliest Trefftz schemes, e.g.\ the \emph{indirect approximation} of \cite[eq.~(35)]{CJZ91}, are of weighted-residual type, even though testing was confined to the boundary of the domain only, see \S\ref{s:OldStyle} below.

\subsubsection{The variational theory of complex rays (VTCR)}
\label{s:VTCR}

The VTCR is a weighted residual Trefftz method introduced in the
1990's by P.~Ladev{\'e}ze and coworkers for problems arising in computational mechanics and later extended to the Helmholtz case in \cite{RLS08}.
Recent surveys are \cite{RLK13,LR14,LBRK12}. 

Several VTCR formulations, slightly different from each other, have been presented.
A general VTCR formulation for the BVP \eqref{eq:BVP} can be written as:
\begin{align}
\text{find}\;u_\Svtcr \in\,& V_p\Th\;\text{s.t.} \quad
\calA_\Svtcr (u_\Svtcr,\vhp)=\ell_\Svtcr (\vhp)\quad
\forall\vhp \in V_p\Th,
\;\text{where}\nonumber\\
\label{eq:VTCRmy}
\calA_\Svtcr (u,v)&:=\im\bigg\{\int_{\Fh^I}\Big(\jmp{u}_N\cdot\mvl{\conj{\nabla_h v}}-\jmp{\nabla_h u}_N\mvl{\conj v}\Big)\di S\\
&
+\int_\GD u\,\conj{\dn v}\di S
+\int_\GR \Big( \frac{C_1}{\ri k\vartheta}(\dn u+\ri k\vartheta u)\conj{\dn v}
+C_2(\dn u+\ri k\vartheta u)\conj v\Big)\di S\bigg\},
\nonumber\\
\ell_\Svtcr (v)&:=\im\bigg\{\int_\GD g_D\conj{\dn v}\di S
+\int_\GR \Big( \frac{C_1}{\ri k\vartheta}g_R\,\conj{\dn v}
+C_2 \,g_R\,\conj v\Big)\di S\bigg\},
\nonumber
\end{align}
where we have reported the formulation with only the imaginary part of the left- and right-hand side, 
following the VTCR convention; however dropping "Im" does not modify the method.

The formulations in \cite[eq.~(21)]{RLK13} and in \cite[eq.~(5)]{LR14}  correspond to the choice of coupling parameters $C_1=1/2$ and $C_2=-1/2$ (up to an overall factor $k$ and using $\re\{-iz\}=\im\{z\}$);
that in \cite[eq.~(6)]{RLSFK12} to $C_1=1/2$ and $C_2=1/2$;
that in \cite[eq.~(4)]{KLR12} to $C_1=1$ and $C_2=0$.
The choice of the coupling parameters does not affect the consistency of the method as all terms in \eqref{eq:VTCRmy} are 
{products} of residuals (internal jumps and boundary conditions) {and traces of test functions}.
In some of the papers cited, using $\im\{a\conj b\}=-\im\{\conj a b\}\,\forall a,b\in\IC$, the conjugation is written on the trial, rather than test, functions in some of the terms, without affecting the formulation.


The VTCR (and similarly the WBM) does not correspond to any of the classical DG schemes listed in \cite{ABC01}.
Indeed, to derive it from the elemental DG equation~\eqref{eq:DGelemental}, one would need to choose  numerical fluxes that, in the terminology of \cite{ABC01}, are neither consistent (they do not equal the fields $\nabla u$ and $u$ when applied to the exact solution $u$ itself) nor conservative (they are not single-valued on the interfaces).


Following \cite[\S2.2]{KLR12}, it is possible to show that if absorption is present then the VTCR is well-posed.
More precisely, provided that $C_1=1$, $C_2=0$, $\re k>0$ and $\im\{k^2\}>0$, the VTCR bilinear form satisfies
$$\calA_\Svtcr (v,v)=-\im\{k^2\}\N{v}_{L^2\OO}^2-\frac{\re k}{|k|^2}\N{\vartheta^{-1/2} \dn v}^2_{L^2(\GR)} \qquad \forall v\in T\Th, $$
thus the VTCR solution is unique in the Trefftz space and coercivity in $L^2\OO$ norm holds (the analogous result for $C_1=-C_2=1/2$ is proved in \cite[Prop.~2]{LR14}).
However, this does not extend to the setting we considered so far,
i.e.\ propagating waves with $k\in\IR$: in this case it can easily be
shown that $\calA_\Svtcr (v,v)=0$ for all $v\in T\Th$ such that $v=0$
on all elements adjacent to the Robin boundary $\GR$ and for any
choice $C_1,C_2\in\IC$, thus well-posedness can not be ensured using a
coercivity argument.
Following \cite[Prop.~2]{LR14}, for $C_1=1/2, C_2=-1/2, k\in\IR$, we have:
$$A_\Svtcr(v,v)=-\frac12\Big(\frac1k\N{\vartheta^{-1/2} \dn u}^2_{L^2(\GR)}+k\N{\vartheta^{1/2} u}^2_{L^2(\GR)}\Big) \qquad \forall v\in T\Th,$$
thus (using Holmgren's theorem \cite[Th.~2.4]{CGLS12}) uniqueness of the solution
of \eqref{eq:VTCRmy} is proved if all mesh elements are adjacent to $\GR$.  
For more general cases, coercivity appears to be too strong an argument. 
We conjecture that discrete inf-sup conditions might be a more viable way for proving well-posedness of the VTCR.

Section 3 of \cite{LR14} considers the application of the VTCR
formulation, corrected with suitable volume terms, with non-Trefftz (piecewise-polynomial) discrete spaces.
This variation is termed \emph{weak Trefftz} and analysed therein.

\subsubsection{The wave based method (WBM)}
\label{s:WBM}

The WBM is a weighted residual Trefftz method, analogous to the VTCR, first introduced in the dissertation of W.~Desmet \cite{DesmetPhD} and later extended to a wide variety of engineering applications, mainly in the realm of vibro-acoustics.
Recent reviews of the state of the art of the research on the WBM can be found in 
\cite{Deckers2014,DesmetChapter2012}.
The discrete space typically used together with the WBM is composed of propagating and evanescent plane waves, as outlined in \S\ref{s:PW}.

The basic variational formulation of the WBM applied to BVP \eqref{eq:BVP}, 
translating \S4.1.4 of \cite{DesmetChapter2012} to our notation and multiplying all terms by $(-\ri k)$, reads
\begin{align*}
\text{find } u_\Swbm\in V_p\Th&\, \text{ s.t.}\quad
\calA_\Swbm(u_\Swbm,\vhp)=\ell_\Swbm(\vhp) \quad
\forall \vhp\in V_p\Th\text{, where}\\
\calA_\Swbm(u,v)&:=\int_{\Fh^I}\bigg(2\jmp{\nabla_h u}_N\mvl{v}
+\frac{\ri k}{Z_{int}}\jmp{u}_N\cdot\jmp{v}_N\bigg)\di S \\
&\qquad+\int_\GR\big(\dn u+\ri k \vartheta u\big)\, v\di S
-\int_\GD u \,\dn v\di S\\
\ell_\Swbm(v)&:=\int_\GR g_R \, v\di S- \int_\GD g_D\,\dn v\di S,
\end{align*}
where $Z_{int}$ is an interior coupling factor.
In some works, a slightly different formulation is used, e.g.\ in \cite[eq.~(81)]{PHV07} different terms are used on the internal interfaces.
We are not aware of any rigorous stability or error analysis of the WBM formulation.



\subsection{Single-element direct and indirect Trefftz methods}\label{s:OldStyle}

%

Most schemes described so far 
were introduced not earlier than mid 1990's, but a lot of research on Trefftz
methods has been carried out since the late 1970's by I.\ Herrera, J.\ Jirousek,
A.P.\ Zieli\'nski, O.C.\ Zienkiewicz and numerous co-workers, mainly for static
elasticity problems.  General reviews of these works are in \cite{Zie97,KK95}; the
Helmholtz case is described in detail in \cite{CJZ91}.  A major difference between
these 
methods and those we described in the previous sections is that in many instances
of the former ones no mesh is introduced on the domain $\Omega$, so that the
unknowns are defined on $\deO$ only.  For this reason, these Trefftz methods more
closely resemble standard boundary element methods rather than finite element
schemes.

There are two main classes of these Trefftz methods: direct and indirect.
(We use the terms ``direct'' and ``indirect'' as in \cite{CJZ91,KK95} and \cite[\S5.1]{PHV07}.)
We describe them for a modification of BVP \eqref{eq:BVP} where we drop the Robin boundary $\GR$ and we consider instead a Neumann boundary portion $\GN$ with boundary condition $\dn u=g_N$.

The \emph{indirect method} is the simplest kind of weighted residual scheme:
\begin{equation}\label{eq:Indirect}
\int_\GD u\,\conj{\dn v}\di S 
-\int_\GN \dn u\,\conj v \di S
= \int_\GD g_D\,\conj{\dn v}\di S 
-\int_\GN g_N \conj v \di S,
\end{equation}
(see \cite[eq.~(35)]{CJZ91} for sound-hard scattering problems in unbounded domains, \cite[eq.~(47)]{PHV07}, \cite[eq.~(16)]{Zie97}, \cite[eq.~(16), (26)]{KK95}).
For Dirichlet exterior problems this is also the method of \cite[\S3]{ADK82}.
In most references the test function is not conjugated.
We note that the indirect method is nothing else than the WBM of \S \ref{s:WBM} posed on a single element, i.e.\ $\calT_h=\{\Omega\}$ and $\Fh^I=\emptyset$.
In the indirect method, the trial functions approximating $u$ are global solutions of the Helmholtz equation on the whole of $\Omega$;
on the other hand the test function $v$ only needs to be defined on $\deO$.
If the Trefftz test and
trial spaces coincide, then the obtained stiffness matrix is symmetric (by
Green's second identity).  If the signs of the terms on $\GN$ are changed, as in
\cite[eq.~(22)]{KK95}, a non-symmetric formulation is obtained.

Subtracting from \eqref{eq:Indirect} the second Green's identity $\int_{\deO}(u\,\conj{\dn v}-\dn u\,\conj v)\di S=0$, which holds for all Helmholtz solutions $u$ and $v$ in $\Omega$, we 
derive the \emph{direct method}:
\begin{equation}\label{eq:Direct}
\int_\GD\dn u\, \conj v\di S-\int_\GN u\,\conj{\dn v}\di S
= \int_\GD g_D\,\conj{\dn v}\di S-\int_\GN g_N\,\conj v\di S, 
\end{equation}
(see \cite[eq.~(42)]{CJZ91}, \cite[eq.~(50)]{PHV07}).  The direct method for the
Dirichlet problem may be viewed as the TDG of \S \ref{s:TDG} with $\alpha=0$ posed
on a single element $K=\Omega$.  Conversely to the indirect method, consistency of
\eqref{eq:Direct} is guaranteed only if the test functions are Helmholtz
solutions in $\Omega$, while the trial functions might be defined (and often are)
on $\deO$ only, for better computational efficiency; the solution is then
evaluated in $\Omega$ with a representation formula in a post-processing step as
for BEMs.
The stiffness matrix arising from the  direct formulation~\eqref{eq:Direct} is the transpose to that of the indirect method~\eqref{eq:Indirect}.
Theorem~6.44 in \cite{Spence2015else} gives sufficient conditions for the well-posedness of the direct method. 
Theorem~7.19 in \cite{SimonSteve2015Chapter} proves that, for well-posed Dirichlet problems with $H^1(\deO)$ data, if the Neumann traces of the trial space coincide with the Dirichlet traces of the test space, then the direct method is well-posed and computes the best approximation of the exact solution in $L^2(\deO)$ norm.
If $\Omega$ is unbounded, the direct and the indirect methods can still be used choosing discrete functions that satisfy Sommerfeld radiation condition; however in \eqref{eq:Direct} the conjugation on the test function must be dropped to preserve consistency.
In this case, if a multipole basis is used, Waterman's \emph{null-field}  method is obtained, see \cite[Ch.~7]{Martin06}, which is a special instance of the \emph{T-matrix} method \cite[\S7.9]{Martin06}.
(Note that \cite{Och95} uses the name \emph{null-field method} for the indirect method with non-conjugated test functions, and \emph{Cremer equations} for the same with conjugated test functions.)

%

For a special choice of Trefftz test functions $v$ indexed by a complex parameter (see the last paragraph of \S\ref{s:PW}), method \eqref{eq:Direct} is called ``\emph{global relation}'' and is the variational formulation at the heart of the \emph{Fokas transform method}, see \cite[eq.~(2)]{DaFo14},
\cite[eq.~(6.142--143)]{Spence2015else} or
\cite[eq.~(7.156)]{SimonSteve2015Chapter}.  
In this context, this formulation 
is typically discretised using piecewise-polynomial (on $\deO$) trial functions, even
though Trefftz functions may be used as well.

\subsection{Non-Trefftz methods with oscillatory basis functions}
\label{s:PUM_et_al}


The main reason for the success of Trefftz methods in the context of time-harmonic
wave problems is that the oscillatory basis functions may offer much better
approximation properties than piecewise polynomials used in standard FEMs.  On the
other hand, similar approximation can also 
be achieved if the discrete functions
are not exact local solution of the PDE to be discretised, but are are only
``approximate solutions''.  If basis functions of this kind are used, the Trefftz
formulations described in the previous sections cannot be 
employed as they stand,
because the residual in the elements will not vanish any more and consistency
will fail.

Approximate Trefftz functions are especially attractive for problems
with smooth\-ly varying material parameters, where no analytic Trefftz function
might be known.  Trefftz formulations, possibly with additional volume terms, can
be used with basis functions that are solutions of the equation only up to a
certain order; see \cite{ImbertGeradDespres2014,TKF14,BeP12}, where this
idea is pursued for DG, UWVF and DEM formulations.

In the following we briefly discuss a few methods that have been proposed employing oscillatory and $k$-dependent basis functions that are not Trefftz.

A very well-known scheme of this kind is the \emph{partition of unity method} (PUM or PUFEM), introduced by I.\ Babu\v{s}ka and J.M.\ Melenk in the mid 1990's, see e.g.~\cite{BAM96}.
The PUM combines the approximation properties of Trefftz functions with the standard variational formulation of the problem, e.g.\ for the BVP \eqref{eq:BVP} with $\GD=\emptyset$
\begin{align}\label{eq:StandardVF}
\int_{\Omega} \big(\nabla_h u \cdot \conj{\nabla_h v}-k^2 u\,\conj v\big)\di V
+\int_\GR \ri k \vartheta u\,\conj v \di S 
=\int_\GR g_R\,\conj v \di S   \qquad \forall v\in H^1\OO.
\end{align}
This requires the use of $H^1\OO$-conforming trial and test functions, thus continuity on interfaces needs to be enforced strongly, which is not viable in Trefftz spaces.
The PUM uses as basis a set of Trefftz functions multiplied to a partition of unity defined on a FEM mesh, e.g.\ piecewise linear/multilinear polynomial FEMs on simplicial/tensor elements.
Theorem~2.1 in \cite{BAM96} ensures that the trial space obtained enjoys the same approximation properties of the Trefftz space employed.
If a $p$-dimensional local Trefftz space is used in each element, together with a piecewise linear/multilinear partition of unity, the total number of degrees of freedom used equals $p$ times the number of mesh vertices, while for a similar Trefftz method on the same mesh (providing comparable accuracy) it would equal $p$ times the number of mesh elements; this means that on tensor meshes almost the same number of DOFs would be employed by the two methods, while on triangles and tetrahedra a saving of a factor up to two or six, respectively, can be achieved by the PUM.
A shortcoming of the PUM is that the formulation \eqref{eq:StandardVF} is not sign-definite and its well-posedness requires a scale resolution condition, while this is not needed for some Trefftz schemes such as the TDG/UWVF presented
in \S\ref{s:TDG} and \S\ref{s:UWVF}.  Differently from Trefftz schemes, the
implementation of the PUM requires the computation of volume integrals; moreover,
the numerical integration of the PUM basis functions may be more expensive than
that of genuine Trefftz functions, see \S\ref{ss:assemble}. 


The PUM for the Helmholtz and other frequency-domain equations was further developed 
by R.J.~Astley, P.~Bettes, A.~El Kacimi, O.~Laghrouche, M.S.~Mohamed, E.~Perrey-Debain, J.~Trevelyan and collaborators, see e.g.\ \cite{LBA02,PLB04}.
When a PUM and a standard FEM discrete spaces are combined, e.g.\ using formulation \eqref{eq:StandardVF}, the method obtained is 
termed \emph{generalised finite element method} (GFEM); e.g.\ \cite{SBH06} employs high-order tensor-product polynomials summed to products of plane waves and bilinear functions.  
In problems with discontinuous wavenumber $k$, the PUM can be applied by coupling
the homogeneous regions by means of Lagrange multipliers as in \cite{LBPT05}; this
is not necessary as formulation \eqref{eq:StandardVF} holds 
on the whole domain,
but enhance the accuracy as in each subdomain only basis functions
oscillating with the correct local wavelength are used.
In \cite{HePfSt08} and related papers, the \emph{trigonometric finite wave elements} (TFWE) is described: the PUM is used with special basis functions adapted to waveguides, lasers and geometries with a single dominant wave propagation direction.
The \emph{finite ray element method} of \cite{MaMa97} consists in the use of a PUM basis in a \emph{first order system of least squares} (FOSLS) formulation; as the unknown is constituted by both $u$ and its gradient, more unknowns are needed but the system matrix is Hermitian.
Finally, in the \emph{hybrid numerical asymptotic method} of \cite{GIK01}, the
PUM space is constructed by multiplying nodal finite elements to oscillating
functions whose phases are derived from geometrical optics (GO) or geometrical
theory of diffraction (GTD), e.g.\ by solving the eikonal equation, \emph{cf.}
\S\ref{ss:adapt}.

Plane wave bases have been combined in~\cite{PPR2015} with the
  \emph{virtual element method} (VEM) framework~\cite{BBCMMR2013}, in order
  to design a high-order, conforming method for the Helmholtz problem, in the spirit of the PUM, but allowing for
  general polytopic meshes.
The main ingredients of the resulting PW-VEM are \emph{(i)} a low
frequency space made of low order
VEM functions, which do not need to be explicitly computed in the element interiors, \emph{(ii)} a proper
local projection operator onto a high-frequency space made of plane
waves, and \emph{(iii)} an approximate stabilisation term. The implementation of the PW-VEM
does not require computation of volume integrals, and no quadrature formulas are required for the
assembly of the stiffness matrix, for meshes with flat interelement boundaries.

The \emph{hybridizable DG} method of \cite{NPRC15} employs two discontinuous discrete spaces (one scalar and one vector) and a space of Lagrange multipliers on the mesh interfaces.
Though Trefftz spaces might be used with this formulation, the authors consider basis functions constructed as products of polynomials and geometrical optics-based oscillating functions, similar to those in \cite{GIK01} but discontinuous.

% ...

A Trefftz approach has been proposed in the context of finite difference schemes in the \emph{flexible local approximation method} (FLAME) by I.\ Tsukerman, 
see e.g.\ the comprehensive review \cite{Tsuk06}. 
In the FLAME, the Taylor expansion of the solution to be approximated used to define classical finite difference schemes is substituted by an expansion in a series of Trefftz basis functions, leading to better accuracy.

Oscillatory basis functions have been successfully used in boundary element methods, in particular for scattering problems, see 
the review on the \emph{hybrid nu\-mer\-i\-cal-asymptotic BEM} (HNA-BEM) \cite{CGLS12}, the \emph{plane-wave basis boundary elements} \cite[\S3]{PLB04} and the
\emph{extended isogeometric boundary element method} (XIBEM) \cite{PTC13}.



\section{Trefftz discrete spaces and approximation}
\label{s:Spaces}

%

Given a Trefftz variational formulation of a BVP, as those in \S\ref{s:Methods}, the definition of a Trefftz finite element method is completed by the choice of a
discrete space 
\[
V_p(\calT_h)=\big\{v\in T(\calT_h):\ 
v_{|K}\in V_{p_K}(K)\big\}\subset T\Th,
\]
where $V_{p_K}(K)$ is a $p_K$-dimensional space of functions $v$ on $K$ such that $\Delta v + k^2 v=0$. 
We describe next the main features of the most common local Trefftz spaces   $V_{p_K}(K)$; 
we do not consider Lagrange multiplier spaces on mesh faces for the methods in \S\ref{s:Farhat}.
The discussion of the conditioning properties of the basis functions described and of the techniques for their numerical integration is postponed to \S\ref{s:Further}.

\subsection{Generalised harmonic polynomials (GHPs)}\label{s:GHP}

Generalised harmonic polynomials 
are smooth Helmholtz solutions that are separable in polar and spherical coordinates in 2D and 3D, respectively, i.e.\ \emph{circular and spherical waves} (also called Fourier--Bessel functions or Fourier basis). 
The local spaces $V_{p_K}(K)$ are defined as follows:
\begin{align*}
&\text{2D:}\quad 
V_{p_K}(K)=\Big\{v:\ v(\bx)=\sum_{\ell=-q_K}^{q_K}\alpha_\ell\, J_\ell(k\abs{\bx-\bx_K})\,\ee^{i\ell\theta},\ \alpha_\ell\in\IC\Big\},
\\
&\text{3D:}\quad 
V_{p_K}(K)=\Big\{v:\ v(\bx)=\sum_{\ell=0}^{q_K}\sum_{m=-\ell}^{\ell}
\alpha_{\ell,m}\,j_\ell(k\abs{\bx-\bx_K})
\,Y_\ell^m\Big(\frac{\bx-\bx_K}{\abs{\bx-\bx_K}}\Big),
\ \alpha_{\ell,m}\in\IC\Big\},
\end{align*}
where $\bx_K\in K$ (e.g.\ is the mass centre of $K$), 
$\theta$ is the angle of $\bx$ in the local polar coordinate system centred at $\bx_K$, $J_\ell$ is the Bessel function of the first kind and order $\ell$, $\{Y_\ell^m\}_{m=-\ell}^{\ell}$ is a basis of spherical harmonics of order $\ell$ (see e.g.~\cite[eq. (B.30)]{AndreaPhD}), and $j_\ell$ is the spherical Bessel function defined by $j_\ell(z)=\sqrt{\frac{\pi}{2z}}\,J_{\ell+\frac{1}{2}}(z)$.
The space dimension $p_K$ is given by $p_K=2q_K+1$ in 2D and by $p_K=(q_K+1)^2$ in 3D.
%
%
We call $q_K$, the maximal index of the (spherical) Bessel functions used, the ``degree'' of the GHP space, as it plays the same role of the polynomial degree in the approximation theory.
A particular feature of GHP spaces is that they are hierarchical. 


The name ``generalised harmonic polynomials'' was coined in \cite{MEL95} and comes from the fact that they are images
of {\em harmonic polynomials} under the operator that maps harmonic functions into Helmholtz solutions, in the framework of Vekua--Bergman's theory~\cite{VEK67,Bergman69} (see also~\cite{VekuaTh,HEN57}). 
The same theory allows to transfer approximation results for harmonic functions by spaces of harmonic polynomials into results on the approximation of Helmholtz solutions by GHPs. 
The density of GHPs in a space of Helmholtz solutions was proved in~\cite[Th.~4.8]{HEN57} and \cite[\S22.8]{VEK67}. 
Approximation estimates in two dimensions were first proved in \cite[Th.~6.2]{EIS74} (in $L^\infty$ norm) and in~\cite{MEL95} (in Sobolev norms), and later sharpened and extended to three dimensions in~\cite{PWapprox}.
We summarise here the estimates of \cite[Th.~3.2]{PWapprox}.


Let $D\in\IR^n$, $n=2,3$, be a bounded, open set with Lipschitz boundary and diameter $h_D$, containing  $B_{\rho h_D}(\bx_D)$ (the ball centred at some $\bx_D\in D$ and with radius $\rho h_D$), and star-shaped with respect to $B_{\rho_0 h_D}(\bx_D)$, where $0<\rho_0\le\rho\le 1/2$. 
Assume that $u\in H^{s+1}(D)$, $s\in\IN$, satisfies $\Delta u+k^2 u=0$ in $D$ and define the $k$-weighted Sobolev norm $\N{u}_{j,k,D}:=(\sum_{m=0}^j k^{2(j-m)}\abs{u}_{m,D}^2)^{1/2}$, $j\in\IN$, where $\abs{\cdot}_{m,D}$ is the Sobolev seminorm of order $m$ on $D$.
\begin{itemize}
\item[{\em i)}] 
If $n=2$ and $D$ satisfies the exterior cone condition with angle $\lambda_D\pi$ \cite[Def.~3.1]{PWapprox}
($\lambda_D=1$ if $D$ is convex), then for every $L\ge s$ there exists a GHP 
$Q_L$ of degree 
at most $L$ such that, for every $j\le s+1$, it holds
\begin{align*}\phantom{.}\hspace{-5.5mm}
\N{u-Q_L}_{j,k,D}\le
C& \big(1+(h_Dk)^{j+6}\big)\ee^{\frac{3}{4}(1-\rho)h_Dk}
\bigg(\Big(\frac{\log(L+2)}{L+2}\Big)^{\lambda_D}h_D\bigg)^{s+1-j}\N{u}_{s+1,k,D},
\end{align*}
where the constant $C>0$ depends only on the shape of $D$, $j$ and
$s$, but is independent of $h_D$, $k$, $L$ and $u$.

\item[{\em ii)}] If $n=3$, there exists a constant $\lambda_D>0$ depending only on the shape of $D$, such that for every $L\ge\max\{s, 2^{1/\lambda_D}\}$ there exists a 
GHP 
$Q_L$ of degree at most $L$ such that, for every $j\le s+1$, it holds
\[
\N{u-Q_L}_{j,k,D}\le
C\big(1+(h_Dk)^{j+6}\big)\ee^{\frac{3}{4}(1-\rho)h_Dk}
L^{-\lambda_D(s+1-j)}
h_D^{s+1-j}\N{u}_{s+1,k,D},
\]
where the constant $C>0$ depends only on the shape of $D$, $j$ and $s$, but is independent of $h_D$, $k$, $L$ and $u$.
\end{itemize}
The main difference between the two results is that the positive shape-dependent parameter $\lambda_D$ entering the exponent of $L$ (thus the $p$-convergence order) is explicitly known in 2D 
(it depends on the largest non-convex corner of $D$) but not in 3D.

Exponential convergence of the GHP approximation of Helmholtz solutions that possess analytic extension outside $D$ were proved in \cite[Prop.~3.3.3]{AndreaPhD} and improved in 2D in~\cite{HMP13}, based upon the corresponding result for harmonic functions of~\cite{HMPS14}. 
Roughly speaking, the error is bounded by a negative exponential of the form 
$C\exp(-bL)\sim C\exp(-bp_{{D}}
^{1/(n-1)})$, while classical bounds for polynomials achieve at most $C\exp(-bp_{{D}}
^{1/n})$, since the dimension $p_D$ of the GHP space of order $L$ is ${\cal O}(L^{n-1})$, while the dimension $p_D$ of the polynomial space of degree $L$ is ${\cal O}(L^{n})$.
Thus, Trefftz methods based on GHPs (and similarly on PWs) can achieve better asymptotic order than standard schemes; however the value of the positive coefficients $b,C$ and their dependence on the BVP and discretisation are not entirely clear.

Approximation estimates in the (discontinuous) spaces $V_p(\calT_h)$ immediately follow from the local approximation estimates with $\overline{D}=K$, for all $K\in\calT_h$.
In case of ($H^1$-conforming) partition of unity spaces enriched with GHPs, global estimates follow from combining the local estimates with~\cite[Th.~2.1]{BAM96}.

GHPs have been proposed in numerous Trefftz formulations:
LS \cite{STO98a,MOW99},
UWVF \cite{LHM13},
VTCR \cite{KLR12},
hybrid-Trefftz \cite[eq.~(62)]{Qin05},
direct and indirect 
single-element schemes \cite{CJZ91,Zie97},
HELS \cite{Wu15},
MPS \cite{FHM67,BeTr05}.

\subsection{Plane waves (PWs)}\label{s:PW}
Plane waves 
probably constitute the most common choice of Trefftz basis functions. 
In this case, the local space $V_{p_K}(K)$ is defined by
\begin{gather}
  \label{s3:pw}
  V_{p_K}(K)=\Big\{v:\ v(\bx)=\sum_{\ell=1}^{p_K}\alpha_\ell\,
  \ee^{\ri k\bd_\ell\cdot(\bx-\bx_K)},\ \alpha_\ell\in\IC\Big\},
\end{gather} where 
$\{\bd_\ell\}_{\ell=1}^{p_K}\subset\IR^n$, $|\bd_\ell|=1$, are distinct propagation directions. 
To obtain isotropic approximations, in 2D, uniformly-spaced directions on the unit circle can be chosen (i.e.\ $\bd_\ell=(\cos (2\pi\ell/p_K), \sin (2\pi\ell/p_K))$); in 3D, \cite{SLW04} and \cite{PTC14} provide directions that are ``almost equally spaced'' (see \cite[\S3.4]{AV05} for a simpler version).
In these cases, the PW spaces are not hierarchical.
However, one of the potential benefits of PW approximations is the possibility to depart from the isotropic case and to adapt the basis propagation directions to the specific BVP at hand and to different elements, either a priori or a posteriori, see \S\ref{ss:adapt}.


The linear independence of arbitrary sets of plane waves (and of their traces) is proved in \cite{AV05,SimonSteve2015Chapter}.
PW bases whose linear independence does not degenerate for small values of $kh_K$ were introduced in \cite[\S3.1]{GHP09} in 2D and in \cite[\S4.1]{PWapprox} in 3D (see also \cite[\S3.4.1]{AndreaPhD}) for analysis purposes.
These stable PW bases converge to GHP bases in the low-frequency limit \cite[p.~815]{PWapprox}.
The existence of these stable bases, which is instrumental to the derivation of approximation estimates for Helmholtz solutions in PW spaces in~\cite{PWapprox}, is guaranteed, provided that the set of directions $\{\bd_\ell\}_{\ell=1}^{p_K}$ constitutes a fundamental system for certain harmonic polynomials. 
In 2D, any choice of $p_K=2q_K+1$ distinct directions, $q_K$ being the maximal degree of the considered harmonic polynomials, guarantees this property. 
In 3D, sufficient conditions on $p_K=(q_K+1)^2$ directions are stated in \cite[Lemma~4.2]{PWapprox}.

Approximation estimates in PW spaces can be derived from similar bounds for GHPs such as those in \S\ref{s:GHP}.
In \cite[Ch.~8]{MEL95}, GHPs were approximated by PWs by approximating their smooth Herglotz kernel with delta functions, leading to $p$-estimates 
in 2D, while in \cite{PWapprox} the Jacobi--Anger expansion was used to link PWs and GHPs in 2D and 3D.
Theorems 5.2 and 5.3 of \cite{PWapprox} (see also \cite[\S3.5]{AndreaPhD}) show that Helmholtz solutions of given Sobolev regularity 
can be approximated in PW spaces with $hp$-estimates similar to those shown in \S\ref{s:GHP} for GHPs.
For PWs, these estimates hold with $L=q_K$, so that $q_K$ plays the role of a ``degree'' for the considered PW space.
As mentioned, for these bounds to hold in 3D, the PW directions have to satisfy some further conditions.
A different derivation of $h$-approximation estimates 
based on a Taylor argument can be found in \cite[Th.~3.7]{CED98}.
In \cite{Per06}, the PW approximation of Helmholtz solutions on the unit disc is analysed in detail, together with the conditioning of different linear systems used for its computation (least squares and collocation for a Dirichlet problem on the disc) and the implications on the accuracy of the approximation computed in finite-precision arithmetic.
We refer again to~\cite[\S5.2]{HMP13} for the exponential convergence in 2D of PW approximations of analytic Helmholtz solutions (see also \cite[Rem.~3.5.8]{AndreaPhD} which holds in 2D and 3D).

Similar to PWs are {the \em evanescent waves}: the basis elements have the same expression $v(\bx)=\ee^{\ri k\bd\cdot\bx}$ but with a more general $\bd\in\IC^n$, $\bd\cdot\bd=1$.
If $\bd=\bd_R+\ri\bd_I$, with $\bd_R,\bd_I\in\IR^n$, then $v$ oscillates in the direction $\bd_R$ (with wavenumber $k|\bd_R|\ge k$) and decays exponentially in the orthogonal direction $\bd_I$ (i.e.\ $|v(\bx)|=\ee^{-k\bd_I\cdot\bx}$). 
Evanescent waves are used in combination with plane waves to approximate interface problems in the DEM~\cite{TZF08} and the UWVF~\cite{LHM13}, and to represent outgoing waves in a 2D unbounded half-strip of the form $\{a<x<b, y>c\}$ in \cite{ZMZ13,SimonSteve2015Chapter}.

A special combination of propagative and evanescent waves is typically used in the WBM.
We describe a 2D version of this space as in \cite[eq.~(14)--(21)]{Deckers2014} (see \cite[\S4.1]{DesmetChapter2012} for 3D).
This 
space is not invariant under rotation but depends on the choice of the Cartesian axes.
For a mesh element $K$, we fix a truncation parameter $N>0$ (typically $1\le N\le 6$) and define $L_x:=\sup_{(x_1,y_1),(x_2,y_2)\in K}|x_1-x_2|$ and $L_y:=\sup_{(x_1,y_1),(x_2,y_2)\in K}|y_1-y_2|$ as the edge lengths of the smallest rectangle containing $K$ and aligned to the Cartesian axes.
Two sets of basis functions are used:
\begin{align*}
&\cos(k_{xj}x)\,\ee^{\pm\ri \sqrt{k^2-k_{xj}^2}\;y}, 
\qquad k_{xj}:=\frac{j\pi}{L_x^K},\quad j=0,\ldots,\lfloor NkL_x^K/\pi\rfloor,\\
&\ee^{\pm\ri \sqrt{k^2-k_{yj}^2}\;x}\,\cos(k_{yj}y),
\qquad k_{yj}:=\frac{j\pi}{L_y^K},\quad j=0,\ldots,\lfloor NkL_y^K/\pi\rfloor,
\end{align*}
for a total dimension $p_K=4+2(\lfloor NkL_x/\pi\rfloor+\lfloor NkL_y/\pi\rfloor)$.
Each basis function is half the sum of two plane (or evanescent) waves, symmetric to one another with respect to the $x$ or $y$ axis: e.g.\ 
$\cos(k_{xj}x)\exp(\ri \sqrt{k^2-k_{xj}^2}y)=\frac12(\ee^{\ri k \bd_{xj}^+\cdot\bx}+\ee^{\ri k \bd_{xj}^-\cdot\bx})$,
with
$\bd_{xj}^\pm:=(\pm k_{xj}/k,\sqrt{1-(k_{xj}/k)^2})$.
A maximum of $4+2(\lfloor kL_x/\pi\rfloor+\lfloor kL_y/\pi\rfloor)$ basis functions are propagative PWs, this number designed to keep the conditioning under control.
If $N>1$, then roughly a fraction $(N-1)/N$ of the total basis functions are evanescent waves decaying in a direction parallel to one of the Cartesian axes.
Refinement is obtained by increasing $N$: 
for $N\le 1$ only propagative waves are present, for higher values evanescent waves are introduced.

In 2D, both evanescent and plane waves may be written as 
$\exp\{\frac k2(\ri (\nu+1/\nu)x+(\nu-1/\nu)y\}=\exp\{\ri k(x\sin\theta+y\cos\theta\}$, parametrised by $\nu\in\IC$ or $\theta\in\IC$ with $\nu=\ee^{\ri\theta}$; these waves constitute the test space (but usually not the trial) for the Fokas method in \cite{Spence2015else,DaFo14} and \cite[\S7.3.4]{SimonSteve2015Chapter} (see also \S\ref{s:OldStyle}).

\subsection{Fundamental solutions and multipoles}\label{s:Hankel}

Fundamental solutions and multipoles are Helmholtz solution in the complement of a point and satisfy Sommerfeld radiation condition 
($\lim_{r\to\infty}r^{\frac{n-1}2}(\der ur-\ri ku)=0$, where $r=|\bx|$).
They are particularly useful to define Trefftz spaces on unbounded elements, e.g.\ for scattering problems.

If the local spaces are spanned by fundamental solutions, simple sources are located at distinct poles $\bx_\ell$ in the complement of $K$:
\begin{align*}
&2D:\quad
V_{p_K}(K)=\Big\{v:\ v(\bx)=\sum_{\ell=1}^{p_K}\alpha_\ell\,
H^{(1)}_0(k\abs{\bx-\bx_\ell}),\ \alpha_\ell\in\IC\Big\},\\
&3D:\quad 
V_{p_K}(K)=\Big\{v:\ v(\bx)=\sum_{\ell=1}^{p_K}\alpha_\ell\, 
\frac{\ee^{-\ri k\abs{\bx-\bx_\ell}}}{\abs{\bx-\bx_\ell}},\ \alpha_\ell\in\IC\Big\},
\end{align*}
where $H^{(1)}_0$ is the Hankel function of the first kind and of order $0$.
Different a priori or a posteriori strategies are used to fix the location of the poles, see \S\ref{s:MFS} and the references cited therein.
As the distance of the points $\bx_\ell$ from $K$ increases, these basis functions approach plane waves, so they permit flexibility not only in the choice of the propagation directions but also in the wavefront curvature. 

Apart from the MFS and its modifications (see \S\ref{s:MFS} and
\cite{FKM03,Och95,BAB08,BAB10,ZMZ14,AV05}), spaces of fundamental solutions
have been used in connection to the UWVF (see \cite{HowarthPhD}, where ray-tracing
is used to determine the poles, and \cite{HowarthChildsMoiola2014}).

Theorem~6 of \cite{Smy09} ensures that Helmholtz solutions in $K$ can be approximated in H\"older norms by 
fundamental solutions centred at any ``embracing boundary'' in 2D and 3D, under weak assumptions on the regularity of~$\deK$.
We are not aware of any result providing orders of convergence.

An alternative approach consists in choosing local spaces generated by multipole expansions, where multiple sources with increasing order are located at a single pole 
$\bx_0$ (or only at few poles):
\begin{align*}
&2D:\;
V_{p_K}(K)=\Big\{v:\ v(\bx)=\sum_{\ell=-q_K}^{q_K}\alpha_\ell\,
H^{(1)}_\ell(k\abs{\bx- \bx_0
})\,\ee^{i\ell\theta},\ \alpha_\ell\in\IC\Big\},
\\
&3D:\;
V_{p_K}(K)=\Big\{v:\ v(\bx)=\sum_{\ell=0}^{q_K}\sum_{m=-\ell}^\ell
\alpha_{\ell,m}\,h^{(1)}_\ell(k\abs{\bx-\bx_0})
\,Y_\ell^m\Big(\frac{\bx- \bx_0}{\abs{\bx-\bx_0}}\Big),\ \alpha_{\ell m}\in\IC\Big\},
\end{align*}
where $H^{(1)}_\ell$ ($h^{(1)}_\ell$) are Hankel functions (spherical Hankel functions, respectively) of the first kind and order $\ell$.
As for the GHPs in \S\ref{s:GHP},  $\theta$ is the angle of $\bx$ in the local coordinate system centred at $\bx_0$, which is located in the complement of $K$, 
and the space dimension is $p_K=2q_K+1$ in 2D and $p_K=(q_K+1)^2$ in 3D.
According to \cite[Rem. 2.2]{BAB10}, fundamental solutions lead to more stable methods than multipoles.

Multipole spaces have been used in connection to 
LS schemes \cite{STO98a,MEH02},
WBM \cite[eq.~(23)]{Deckers2014}, \cite[\S4.1.2]{DesmetChapter2012},  
hybrid-Trefftz \cite[eq.~(63)]{Qin05}, 
HELS \cite{Wu15},
source simulation techniques \cite{Och95},
null-field \cite{Martin06} and
single-element schemes \cite{CJZ91,Zie97,ADK82}.
In \cite{HBB99} and related papers, some 2D multipoles with suitably chosen index $\ell$ (not necessarily integer) are used on infinite sectors, in such a way to ensure continuity of discrete functions across rays; this might be more efficient than full multipole spaces 
for solutions with a preferred propagation direction.


\subsection{Other basis functions}\label{s:OTHERS}

Other discrete Trefftz spaces have been proposed in literature for use with the
  various approaches covered in \S\ref{s:Methods}.

In 2D, {\em corner waves} such as $J_{\ell/\alpha}(k|\bx|)\sin(\ell\theta/\alpha)$, with $\ell\in\IN$ and $0<\alpha<2$, capture the behaviour of Helmholtz solutions near a domain corner of angle $\pi\alpha$.
They have been used e.g.\ in the WBM \cite{DBVVD12},    
in LS methods~\cite{BAB10,STO98a,ZMZ13} and in the MPS \cite{FHM67,BeTr05}.    
In \cite{ZMZ14}, they are used with $\alpha=2$ on tips of 1D screens to represent the strong singularities of the solution in a non-Lipschitz domain.
Theorem 6.3 of \cite{EIS74} uses Vekua--Bergman theory to give orders of convergence for the approximation of singular functions by spaces of corner waves and GHPs (see also \cite[\S5]{BAB10} and references therein).
We are not aware of any use of similar functions in 3D.

%
%

The {\em wave band functions}, introduced in the VTCR context \cite{RLS08}, are Herglotz functions with piecewise-constant kernel, e.g.\ 
$\int_a^b \ee^{\ri k(x\cos\theta+y\sin\theta)}\di \theta$ in 2D.

In the presence of a circular hole, suitable combinations of {\em Hankel and Bessel functions} a priori fulfil homogeneous boundary conditions~\cite[eq.~(13)]{STO98a}.

If the wavenumber varies inside an element, the basis functions described so far do not lead to Trefftz methods.
In case of linearly variable wavenumber, {\em Airy functions} can be used to construct Trefftz spaces \cite{TKF14}.
In \cite{ImbertGeradDespres2014,IG15} \emph{generalised plane waves} in the form $\ee^{P(\bx)}$, for suitable polynomials $P$, are introduced and analysed in a UWVF setting: they solve a perturbed Helmholtz problem and converge with high orders in $h_K$.
Similar  ``almost-Trefftz'' waves are used in~\cite{GDA07} and named \emph{oscillated polynomials}. 
\emph{Modulated plane waves}, i.e.\ products of PWs and polynomials, are the basis functions of the DG method of \cite{BeP11,BeP12}; as they are only ``approximately Trefftz'', volume terms appear in the formulation.

Products of (continuous) low-order polynomials and PWs or GHPs constitute the basis of the PUM \cite{BAM96,PLB04,SBH06,LBPT05,HePfSt08}, while products of polynomials and oscillating functions derived from high-frequency asymptotics are basis elements in \cite{GIK01,NPRC15}.

\section{Further topics}
\label{s:Further}

\subsection{Assembly of linear systems}
\label{ss:assemble}

All the Trefftz finite element methods for \eqref{eq:BVP} discussed in
\S\ref{s:Methods} give rise to dense or sparse linear systems of
equations. Entries of coefficient matrices are obtained by integrating products of
(derivatives of) trial and test functions over bounded $d$-dimensional 
sub-manifolds of $\Omega$, $d<n$. The stable and accurate (approximate) evaluation
of these integrals is a key implementation issue. 

Among all Trefftz approximation spaces and associated bases presented in
\S\ref{s:Spaces}, plane waves (PWs) $\ee^{\ri k\bd\cdot\bx}$ (either propagative with $\bd\in\IR^n$ or evanescent with $\bd\in\IC^n$) are exceptional, because they allow a closed-form evaluation of their integrals over any flat sub-manifold with piecewise
flat/straight boundary. 
For instance, in all variants of PW-based Trefftz methods on polyhedral meshes in 3D,
expressing mesh faces by 2D parametrisations,
we eventually encounter integrals of the form
\begin{gather}
  \label{s4:1}
  \int\nolimits_{F} \exp(\bw\cdot \bx)\di V,\qquad
  F \subset \mathbb{R}^{2}
  \;  \text{a bounded polygon, $\bw\in\IC^2$ constant.}
\end{gather}
Then we can take the cue from \cite[\S2.1]{GAB09} or \cite[\S4]{KLA10} and apply integration
by parts in order to reduce \eqref{s4:1} to integrals over the straight edges
$e_{1},e_{2},\ldots e_{q}$, $q\in\mathbb{N}$ of $F$:
\begin{gather*}
  \int\nolimits_{F} \exp(\bw\cdot \bx)\di V = 
   \frac{1}{{\bw\cdot\bw}}\int\nolimits_{F} \bw\cdot \nabla\exp(\bw\cdot \bx)\di V = 
  \sum_{\ell=1}^{q} \frac{\bw\cdot\bn_{\ell}}{{\bw\cdot\bw}} \int\nolimits_{e_{\ell}} 
  \exp(\bw\cdot\bx)\di  s,
\end{gather*}
where $\bn_{\ell}$ is the exterior normal at $e_{\ell}$. Then, as in
\cite[Ch.~2]{GIT08}, if $e_{\ell}=[\ba,\bb]$, $\ba,\bb\in\mathbb{R}^{2}$, we find,
$\int\nolimits_{e_{\ell}} \exp(\bw\cdot\bx)\di {s} 
= \exp(\bw\cdot\ba)|\bb-\ba|
\psi(\bw\cdot(\bb-\ba))$, where $\psi(z) := (\exp(z)-1)/z$. Of course, a
numerically stable implementation of this function for small arguments is
essential\footnote{A stable algorithm for point evaluations of $\psi$ even for
  arguments close to $0$ is provided by the MATLAB function \texttt{expm1}.}. This
approach can be generalised to yield analytic formulas for computing integrals of
products of PWs times polynomials, see \cite{KLA10,GAB09}, with increased
computational effort, however. 

Approximate evaluation of the integrals becomes inevitable for all choices of Trefftz basis functions {other than} PWs, and even for a PW basis on meshes with curved elements. 
Then Gauss--Legendre numerical quadrature
seems to be the most widely used option. However, the integrands 
may be oscillatory, which delays the onset of (exponential) convergence of the
quadrature error until the number of quadrature points surpasses a threshold
roughly {proportional to} the ratio of the local mesh size and the wavelength. This leads
to higher computational cost per degree of freedom for larger values of $kh_{K}$.
One may think of using special quadrature rules for oscillatory integrals, as
derived, for instance, in \cite{HUO08}. Those avoid an increase in the number of
quadrature points for growing spatial frequency of the oscillations, but
unfortunately 
require precise knowledge of the oscillatory term in the integrand.



\subsection{Adaptive Trefftz methods}
\label{ss:adapt}

Besides classical $h$-, $p$- or $hp$-adaptivity, Trefftz methods offer scope for more sophisticated adaptive strategies consisting in the choice of specific basis functions for different BVPs and in different mesh elements, either a priori or a posteriori.

The main strand of \textbf{a priori adaptive} Trefftz methods falls into the
category of \emph{hybrid numerical-asymptotic} methods. 
High-frequency limit models, such as ray optics or geometric theory of diffraction (GTD), guide the
selection of local Trefftz spaces in the individual cells of a mesh. In a
non-Trefftz PUM framework this idea was pursued in \cite{GIK01}, and within the hybridizable DG method in \cite{NPRC15}, in both cases for 2D acoustic
scattering at a smooth sound-soft object. 
In these works, local phase factors $\bx\mapsto\exp(\ri k S(\bx))$ derived from reflected and diffracted waves multiply standard continuous nodal basis functions, 
in \cite{GIK01}, or local polynomials, in \cite{NPRC15}, thus generating a basis for (local) trial spaces.

The policy of incorporating local directions of rays is particularly attractive for 
PW-based methods,
because PW basis functions naturally encode
a direction of propagation. For problems where excitation is due to an incident
PW and material properties are piecewise constant, ray tracing and related
techniques \cite[\S3.2]{NPRC15} based on geometric optics (specular reflection and
Snell's law of refraction at material interfaces) can provide information about
the local orientation of wave fronts for $k\to\infty$. PWs matching the
found ray directions are then used to build local bases, either exclusively or
augmented by a reduced set of generalised harmonic polynomials (GHPs) or ``equi-spaced'' PWs.

This idea for TDG was first outlined and tested in \cite{BeP11} and further
elaborated and extended in \cite[Ch.~5]{HowarthPhD} (for UWVF). In the latter
work, in an attempt to resolve curved wave fronts and take into account diffracted
waves from corners, also Hankel functions $\bx\mapsto
H^{(1)}_{0}(k|\bx-\by_{\ast}|)$ with $\by_{\ast}$ outside a mesh cell have been
proposed as local basis functions. Approximation of curved wave fronts deduced
from GTD corrections is also attempted in \cite{BeP12}. There the authors move
beyond Trefftz methods and use DG with trial spaces of polynomially modulated
PWs, which are more suitable for approximating propagating circular waves. 

In simple 2D situations with convex smooth or polygonal scatterers and incident
plane wave, overall accuracy seems to benefit substantially from a priori
directional adaptivity. However, if there are more than only a few dominant wave
directions as in the case of more complicated geometries, trapping of waves, dark
zones and shadow boundaries, current directional adaptivity soon meets its
limitations. On the other hand, this strategy appears as the most promising way to achieve
$k$-uniform accuracy with numbers of degrees of freedom that remain $k$-uniformly
bounded or display only moderate growth as $k\to\infty$.  The potential of this
idea has been strikingly demonstrated in the case of BEM for 2D scattering
problems \cite{CGLS12}.

\textbf{A posteriori directional adaptivity} seeks to extract information about
dominant wave directions from intermediate approximations of $u$.  A
refine-and-coarsen strategy is embraced in \cite{BeP11}. In each step of the
adaptive cycle it first computes a PWDG solution $u$ of the scattering problem
based on a relatively large number of local Trefftz basis functions (GHPs and PWs). Subsequently,
by solving local non-linear $L^{2}$-least squares problems, the directions of fewer
PWs are determined so that $u$ can still be well approximated locally. 

A $p$-hierarchical error indicator is studied in \cite{GIT08}. In a step of the
adaptive scheme starting from the approximate solution $u$ a presumably improved
solution $\hat{u}$ is computed using double the number of local PWs. Then
a single local plane wave direction $\bd_{K}$ on a mesh element $K$ is extracted
from the error $e(\bx):=\hat{u}(\bx)-u(\bx)$ through the projection formula
\begin{align*}
  \widetilde{\bd}_{K} := \re \int_{K}\frac{\nabla e(\bx)}{\ri k e(\bx)}\di V,\qquad
  \bd_{K} := \frac{\widetilde{\bd}_{K}}{|\widetilde{\bd}_K|}.
\end{align*}
Detailed numerical experiments are reported in \cite[Ch.~6]{GIT08}. In the
pre-as\-ymp\-tot\-ic regime, when the resolution of the trial spaces is still rather
low, one observes a pronounced gain in accuracy in the case of the adaptive
approach compared to approximation with the same total number of equi-spaced
PWs.

Directional adaptivity for Trefftz methods has also been tried in other flavours.
In the context of least squares methods as discussed in \S\ref{s:LS} an
offset angle for the sets of local equi-spaced PWs is introduced as another
degree of freedom in \cite{ACD13}, aiming to align them with a local dominant wave
direction. 
For the VTCR method presented in \S\ref{s:VTCR},
an error indicator based on local wave energy is used in \cite{RLSFK12} 
to steer angular refinement of local Trefftz spaces.

\textbf{A posteriori mesh adaptivity} is considered in \cite{KMW14}, where classical ``elliptic'' error estimation and mesh refinement strategies are adapted
for the $h$-version of TDG. In a low-frequency setting, the method inherits the
good performance of the underlying adaptive mesh refinement algorithms for polynomial
DG for the Poisson equation. However, there is little hope that this carries over 
to larger wavenumbers~$k$. 
A similar error estimator, aimed at adaptive mesh refinement, has been
  described in \cite[\S3.2]{ADF09} for the DEM/DGM presented in \S\ref{s:Farhat}.


\subsection{Ill-conditioning and solvers}
\label{ss:solvers}

The linear systems of equations spawned by PW-based finite element methods
are highly prone to ill-conditioning, when high resolution trial spaces are used,
see e.g.~\cite[\S5]{HMK02}, \cite[\S4.3]{GAB07}, \cite{GAA07}, and \cite{LBA02}
for a PUM setting. This is largely caused by an inherent instability of the PW
basis on cells, whose size is relatively small compared to the
wavelength. Intuitively, for $|\bx| \ll k^{-1}$, the functions $\bx\mapsto
\ee^{\ri k \bd_{\ell}\cdot\bx}$ from \eqref{s3:pw} are almost constant, hence,
nearly linearly dependent, \emph{cf.} \cite[\S4.2]{LBA02}. The same heuristics
applies, when their density increases; even for cell sizes comparable to the
  wavelength, PWs are hardly distinct when their directions are close,
\emph{cf.} \cite[\S4.3]{LBA02}.

Empirically, for the local PW Galerkin matrix $\mathbf{M}_{K}$ associated
with the $L^{2}$ inner product on a single mesh cell $K$,
we find that its spectral condition number grows like $\sim h_K^{-q}$ for cell size $h_{K}\to 0$,
  where $q>0$ is 
  proportional to the number $p_{K}$ of (approximately
  uniformly spaced) PWs in 2D, and to the square root of $p_K$ in 3D.
  Essentially, $q$ is related to the ``degree'' of the considered set of $p_K$ 
  PWs; see \S\ref{s:PW}.
  Even worse, the condition number soars exponentially in $q$: 
  $\operatorname{cond}_{2}(\mathbf{M}_{K}) \sim \ee^{\alpha q}$ for
  $q\to\infty$ and $\alpha>0$; see Appendix~A.
%
A similar explosion of condition numbers is
observed for the full systems matrices as meshes are refined or more PW
basis functions per element are used.

There is circumstantial evidence that direct sparse elimination can cope fairly
well with the ill-conditioned linear systems arising from UWVF or PUM, see
\cite[\S5.3.3]{GAA07}, \cite{LHM13}. Yet, eventually the instability of the basis
will impact the quality of the solution \cite[\S5.4]{SBH06}. A remedy proposed in
\cite{HMK02} for the UWVF is to limit $p_{K}$ based on monitoring condition
numbers of element matrices. Apparently, this also curbs the condition number of
the global system matrix. Alternatively, there exist different heuristic recipes
for choosing a priori the number of PWs per element to balance accuracy and
conditioning: in 2D, the widely cited \cite[eq.~(14)]{HGA08} suggests
$p_K=\mathrm{round}(kh_K+C(kh_K)^{1/3})$ with $3\le C\le 14$ for the UWVF, while
\cite[\S5.1.1]{LR14} proposes $p_K=\lfloor 2kh_K\rfloor$ for the VTCR.
For the WBM, \cite[\S3.2]{Deckers2014} proposes a rule to balance propagative and evanescent basis functions, see \S\ref{s:PW}.

The most straightforward cure for instability would trade the PW basis of
$V_{p_{K}}(K)$ from \eqref{s3:pw} for a more stable basis, found by local
orthonormalisation as in the case of polynomial FEM, \emph{cf.}
the approach from \cite[\S3.1]{NPRC15}. However, instability may sneak in through
the back door and manifest itself in severe impact of round-off errors during
orthonormalisation and recombination of element matrices. The use of high-precision arithmetic may be advisable, but has never been documented.

For the sake of stability, PWs may be replaced by the generalised harmonics polynomials introduced in \S\ref{s:GHP}. 
In 2D, a scaling of the GHPs has been devised in~\cite{LHM13}, in order to lower the condition number of the resulting UWVF:
\begin{gather*}
  \frac{J_\ell(k\abs{\bx-\bx_K})\,\ee^{i\ell\theta}}{k\sqrt{\abs{J^\prime_{\ell}(kh_K)}^2+\abs{J_{\ell}(kh_K)}^2}}.
\end{gather*}
In \cite{LHM13},  
it is also shown that the conditioning of GHP-based UWVF schemes is better than for methods based on PWs, and that it improves on regular meshes. 
This might be related to the orthogonality of GHPs on balls.

The numerical experiments in \cite[\S3.7]{HowarthPhD} suggest that the use of fundamental solutions as basis functions may considerably reduce the conditioning of UWVF matrices, at the expense of accuracy.
Both accuracy and conditioning increase the further the centres of the fundamental solutions are from the element.

The use of iterative solvers for linear systems generated by 
Trefftz methods entails preconditioning. For PW basis functions, the 
first proposal in \cite[\S2.4]{CED98} for the UWVF was a local preconditioner, equivalent to an orthonormalisation of the PW basis with respect to an $L^{2}$ inner product on the boundary of mesh cells. 
An interesting connection of the local preconditioner with non-overlapping
optimised Schwarz domain decomposition methods was discovered in \cite{BGP14}.
The local preconditioner was used in conjunction with a BiCGStab Krylov subspace
solver in \cite{HMK02} and augmented by a coarse-grid correction in the spirit of
non-overlapping domain decomposition in \cite{HUY13,HuYuan14}. The coarse space is
again spanned by PWs. This is also true for the two-level sub-structuring
preconditioner proposed for DEM/DGM (see \S\ref{s:Farhat}) in \cite{FTT09}.
Two-level, non-overlapping Schwarz domain decomposition preconditioners for
  PWDG (essentially UWVF) have been tested in~\cite{APZ15}; these preconditioners
  seem to be robust with respect to the wavenumber $k$ and the local number
  of PW directions, although they do not seem to be perfectly scalable with
  respect to the number of subdomains. 

\section{Assessment and conclusion}
\label{s:concl}

Faced with a flurry of different Trefftz methods and a wealth of numerical data,
we feel at a loss about making unequivocal statements about the merits of Trefftz
methods, let alone ranking them according to some undisputed criteria. Rigorous
theory is available for LS methods (\S\ref{s:LS}), TDG (\S\ref{s:TDG}),
and PUM (\S\ref{s:PUM_et_al}). 
Combined with approximation results for suitable Trefftz bases, this leads to better asymptotic estimates in terms of orders of convergence in the number of degrees of freedom to what is available for polynomial FEM (e.g.\ \cite{PVersion,HMP13}).
The dependence of crucial constants on the wavenumber $k$ is explicit in several cases, but the orders in $k$ are usually not better than for polynomial methods.
Thus theory fails to provide information about the key issue of
``$k$-robust'' accuracy with ``$k$-independent'' cost. Moreover, numerical
dispersion will also haunt local Trefftz methods in the case of $h$-refinement;
thus they provide no escape from the pollution error. 

We also advise caution when reading numerical experiments, because they may be
tarnished by selection bias, making authors 
subliminally 
pick test cases matching the intended message of an article. Disregarding this, even
``objective'' comparisons are inevitably confined to a few simple model
problems. This is problematic, because different model problems
sometimes seem to support opposite conclusions. 

From our experience, the power of Trefftz methods can best harnessed by
$p$-refinement using approximation by Trefftz functions in regions as large as
possible.  In the presence of singularities we recommend either the use of
  corner basis functions (\S\ref{s:OTHERS}) in 2D, or $hp$-refinement, maybe using standard polynomial approximation
on small elements as in \cite{MSS10}. There is a solid theoretical
foundation, when this is done in the LS, TDG, or PUM framework. The resulting
methods should be able to compete successfully with polynomial FEM even in their
more sophisticated versions tailored to wave propagation problems
\cite{EM11,MES11,FEW11}.

The discussion of adaptive approaches in \S\ref{ss:adapt} hints that some Trefftz
trial spaces have approximation capabilities well beyond the reach of polynomials. 
Directional adaptivity seems to be very promising, but much research will still
be required to convert it into a reliable practical algorithm. The same applies
to iterative solvers and preconditioners for Trefftz schemes, see
\S\ref{ss:solvers}, which might also benefit considerably from the extra
information contained in Trefftz trial spaces. Hence, we believe that many exciting
possibilities offered by the idea of Trefftz approximation still await discovery 
and that the full potential of Trefftz methods is only gradually being realised.


%
%
%



\section*{Appendix A: Condition numbers of plane wave mass matrices}
\label{app:cond}

\begin{figure}[b]
  \centering
  \begin{minipage}[c]{0.5\textwidth}\centering
    \includegraphics[width=0.95\textwidth]{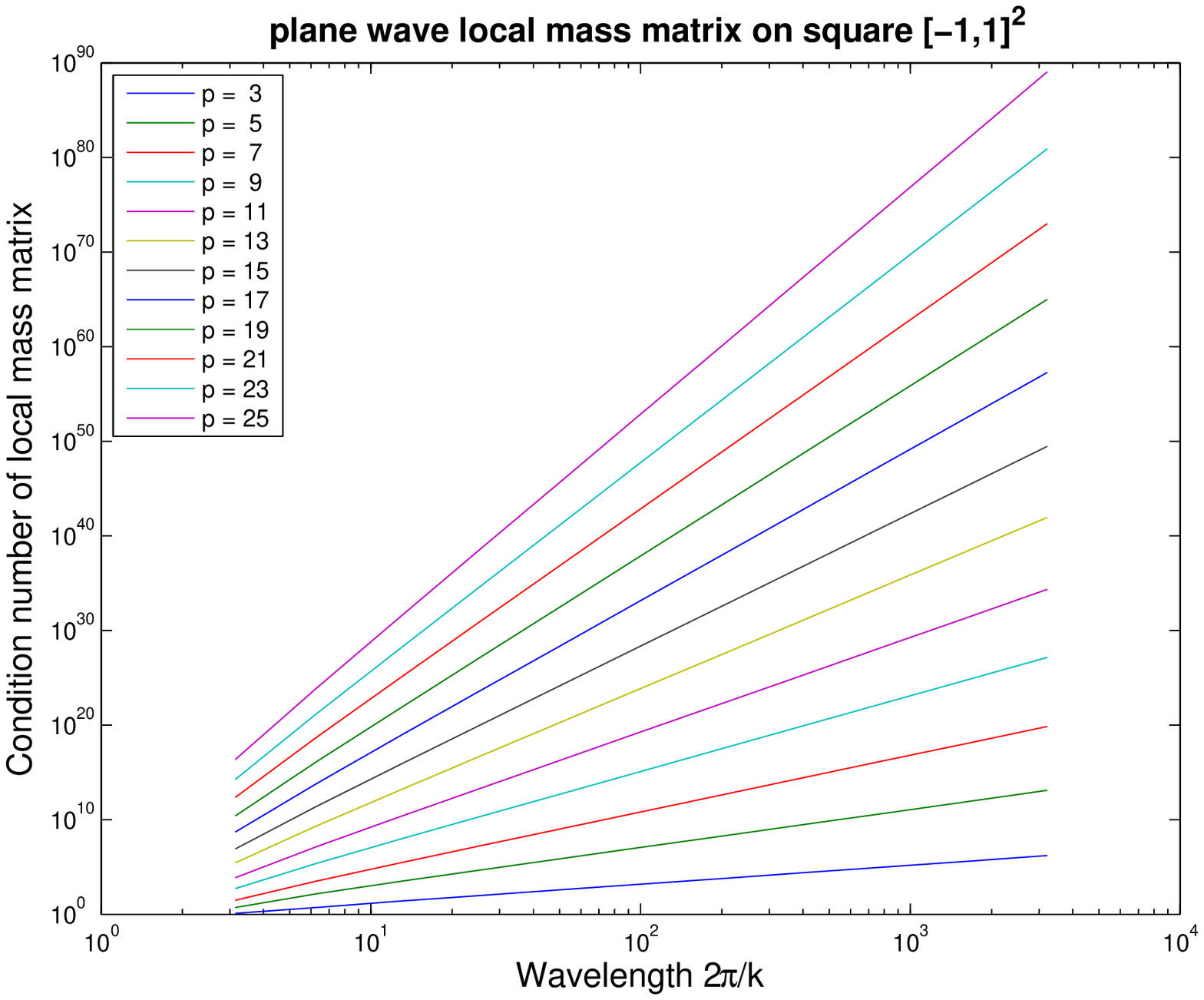}
  \end{minipage}%
  \begin{minipage}[c]{0.5\textwidth}\centering
    \includegraphics[width=0.95\textwidth]{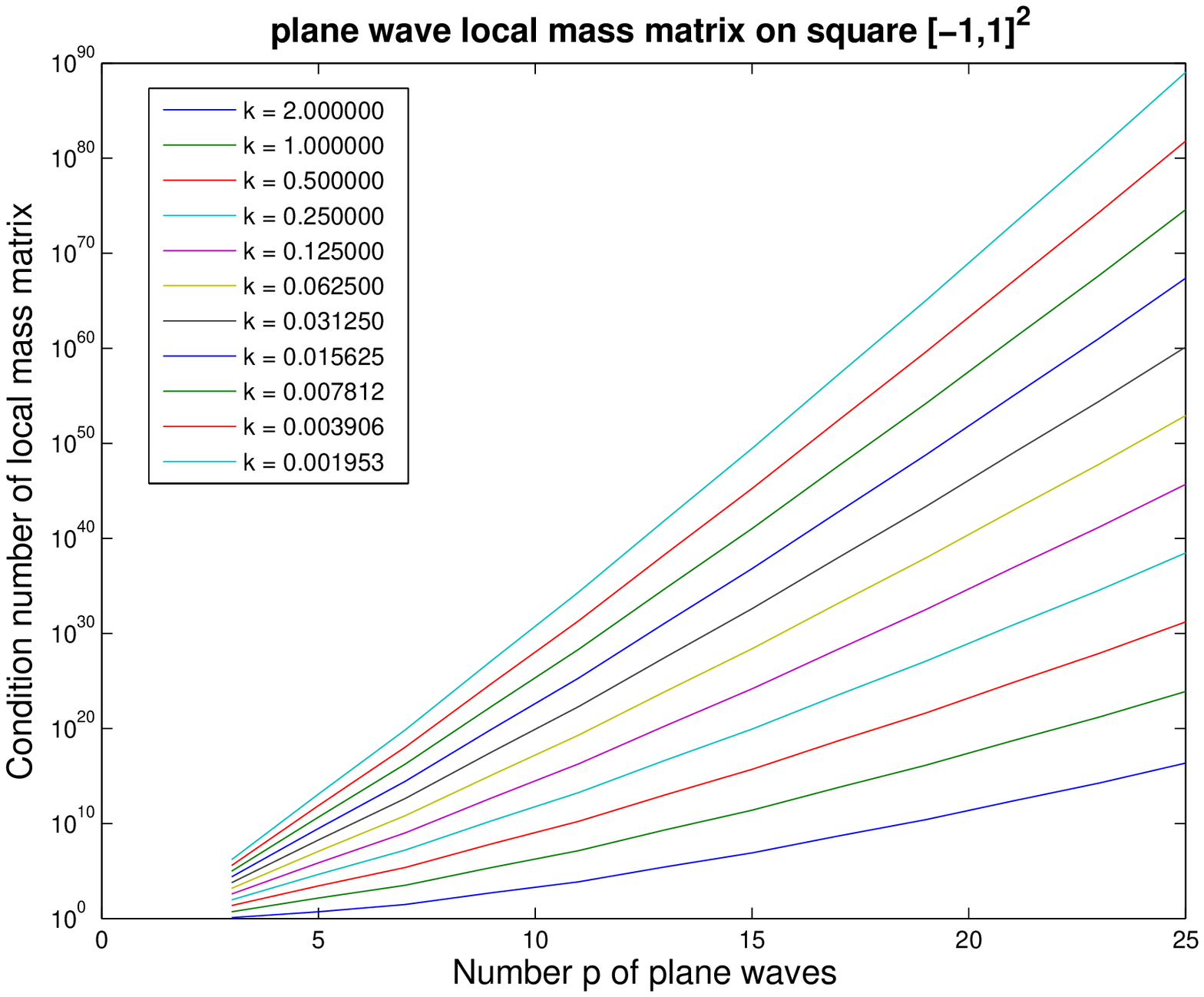}
  \end{minipage}%
  \caption{Condition numbers of element mass matrices on the square $(-1,1)^2$}
  \label{fig:cond2d}
\end{figure}

\begin{figure}[b]
  \centering
  \begin{minipage}[c]{0.5\textwidth}\centering
    \includegraphics[width=0.95\textwidth]{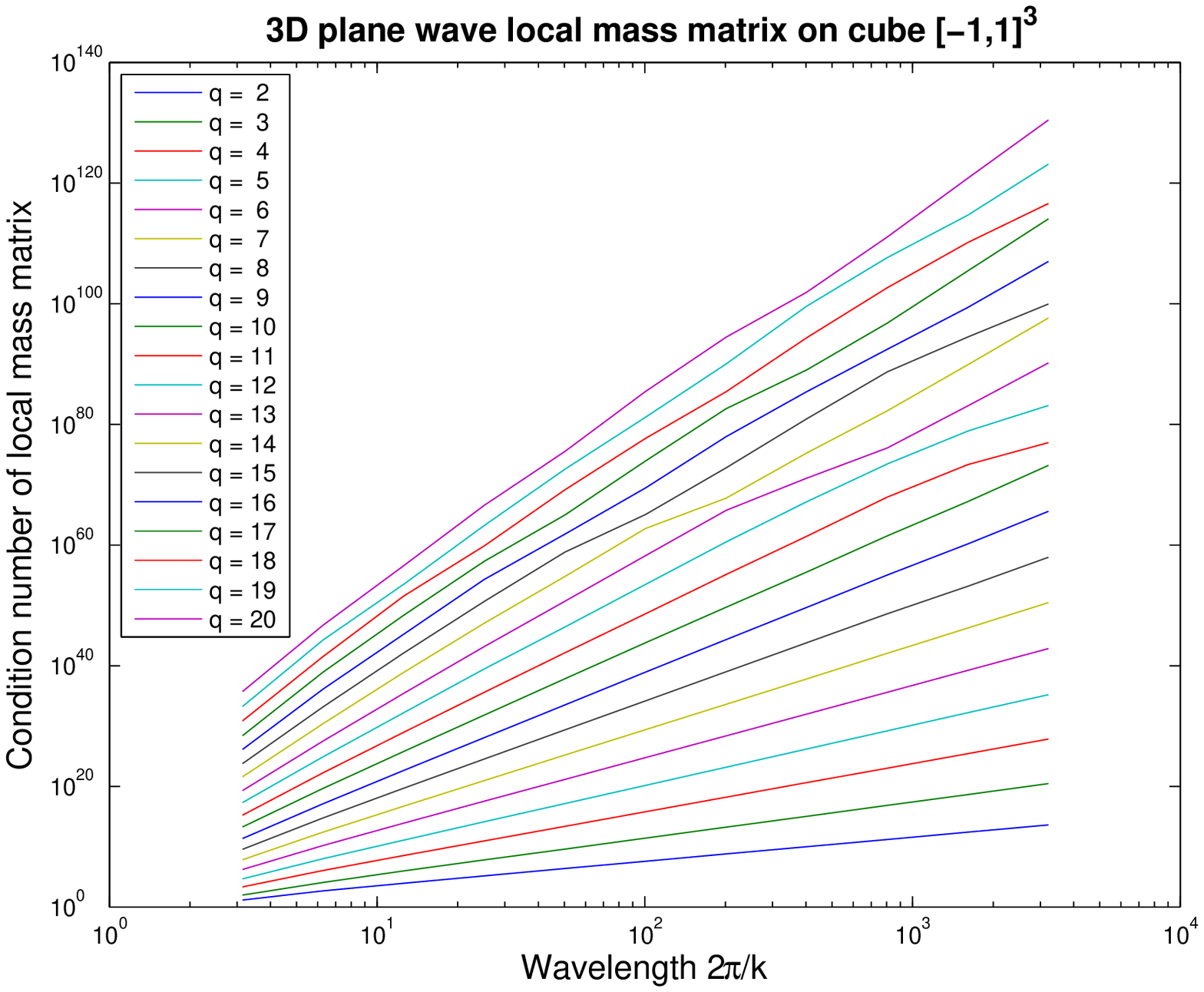}
  \end{minipage}%
  \begin{minipage}[c]{0.5\textwidth}\centering
    \includegraphics[width=0.95\textwidth]{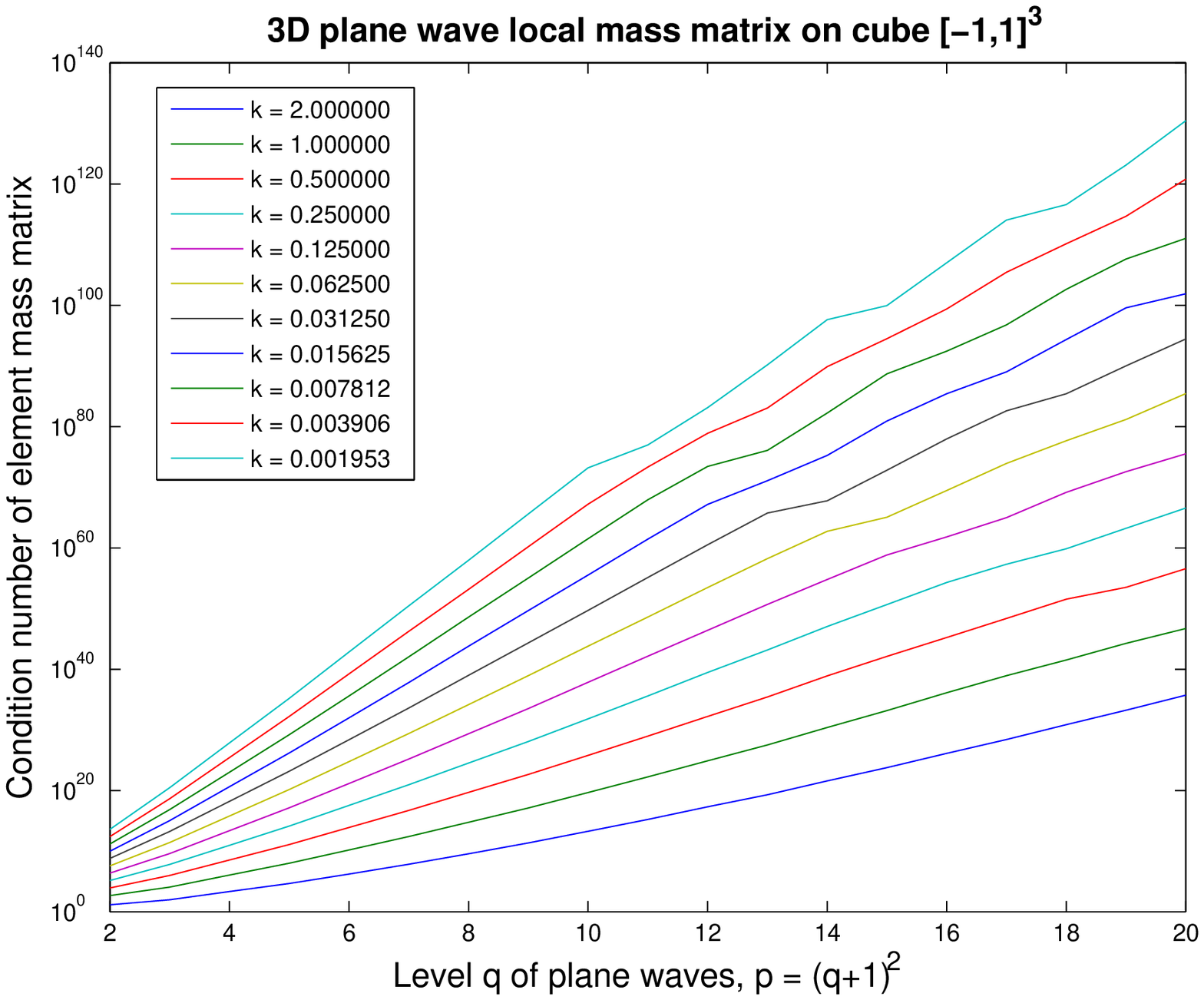}
  \end{minipage}%
  \caption{Condition numbers of element mass matrices on the cube $(-1,1)^3$. }
  \label{fig:cond3d}
\end{figure}

\newcommand{\wn}{k}
Given a wave number $\wn>0$ and $p\in\mathbb{N}$ distinct unit vectors $\bd_{\ell}\in\mathbb{R}^{n}$,
$\ell=1,\ldots,p$, and a domain $K\subset\mathbb{R}^{n}$ with barycentre $\bx_K$, the symmetric positive
definite plane wave element mass matrix $\mathbf{M}_{K}$ on $K$ is defined as
\begin{gather*}
  \mathbf{M}_{K} := \left(\int\nolimits_{K} e^{\ri\wn \bd_{\ell}\cdot(\bx-\bx_K)}\cdot e^{-\ri\wn
    \bd_m\cdot(\bx-\bx_K)}\di V\right)_{\ell,m=1}^{p}.
\end{gather*}
For $n=2$ we computed spectral condition numbers of $\mathbf{M}_{K}$ for equi-spaced directions
$\bd_\ell=(\cos (2\pi\ell/p), \sin (2\pi\ell/p))$, $\ell=0,\ldots,p-1$. For $n=3$
we choose the directions $\bd_{\ell}$ as the ``minimum norm points'' 
according to
I.H.~Sloan and R.S.~Womersley \cite{SLW04,SLW--}. These points are indexed by a level
$q\in\mathbb{N}$ and $p=(q+1)^{2}$.  The spectral condition numbers are plotted in
Figure~\ref{fig:cond2d} for $n=2$, $K=(-1,1)^{2}$, and Figure~\ref{fig:cond3d} for
$n=3$, $K=(-1,1)^{3}$.  They have been computed with MATLAB using the high-precision arithmetic (200 decimal digits) provided by the Advanpix Multi-Precision
Toolbox\footnote{\url{http://www.advanpix.com/}}.


\end{document}